\documentclass{article}
\usepackage{amssymb}
\usepackage{amsmath}
\usepackage{hhline}

\input{tcilatex}
\begin{document}

\begin{center}
{\LARGE Commuting Magic Square Matrices\medskip }

{\large Ronald P. Nordgren\footnote{%
email: nordgren@rice.edu}}

Brown School of Engineering, Rice University{\LARGE \medskip }
\end{center}

\noindent \textbf{Abstract.} We review a known method of compounding two
magic square matrices of order $m$ and $n$\ with the all-ones matrix to form
two magic square matrices of order $mn.$ We show that these compounded
matrices commute. Simple formulas are derived for their Jordan form and
singular value decomposition. We verify that regular (associative) and
pandiagonal commuting magic squares can be constructed by compounding. In a
special case the compounded matrices are similar. Generalization of
compounding to a wider class of commuting magic squares is considered. Three
numerical examples illustrate our theoretical results.

\section{Introduction}

The construction of magic squares by compounding smaller ones arose over
1000 years ago and has attracted interest ever since as reviewed by Pickover 
\cite{PICK}, Chan and Loly \cite{CHAN}, and Rogers, et.\negthinspace\ al 
\cite{ROGE}. In the present article we begin with a formulation of the
compounding construction given by Eggermont \cite{EGGE} and extended in \cite%
{ROGE}. Two basic magic square matrices of order $m$ and $n$\ $\left( m,n\
\geq 3\right) $ are compounded with the all-ones matrix by means of the
tensor (Kronecker) product to form two unnatural magic square matrices of
order $mn$ that are shown to commute. These two matrices are combined to
form two natural magic square matrices of order $mn$ that commute. When $m=n$
and the two basic magic squares are identical, we verify a result in \cite%
{ROGE} that the two compounded magic square matrices are related by a
row/column permutation (shuffle) that we express in matrix form which shows
that they are similar matrices. We verify a known result \cite{EGGE} that
regular (associative) and pandiagonal special properties are preserved by
compounding. Simple formulas are derived for the matrices in the Jordan form
and singular value decomposition of the compounded matrices in terms of
those of the basic matrices . Generalization to a wider class of commuting
magic squares is considered but our formulas for the Jordan form and SVD do
not apply to them. Examples are given of commuting magic square matrices of
orders 9, 12, and 16, all with special properties.

\section{Construction}

\noindent \textbf{Magic Squares. }To begin, let $\mathbf{E}_{n}$ denote the
order-$n$ square matrix with all elements one and let $\mathbf{R}_{n}$
denote the order-$n$ square matrix with ones on the cross diagonal and all
other elements zero. Let $\mathbf{M}_{n}$ be an order-$n$ \emph{magic square 
}matrix whose rows, columns and two main diagonals add to the \emph{%
summation index} $\mu _{n},$ i.e.%
\begin{align}
\mathbf{M}_{n}\mathbf{E}_{n}& =\mathbf{E}_{n}\mathbf{M}_{n}=\mu _{n}\mathbf{E%
}_{n},  \notag \\
\limfunc{tr}\left[ \mathbf{M}_{n}\right] & =\limfunc{tr}\left[ \mathbf{R}_{n}%
\mathbf{M}_{n}\right] =\mu _{n}.  \label{MnMagic}
\end{align}%
We further require that $\mathbf{M}_{n}$ be a \emph{natural} magic square
with elements $0,1,\ldots ,n^{2}-1$, whence%
\begin{equation}
\mu _{n}=\frac{n}{2}(n^{2}-1).  \label{mn}
\end{equation}%
From two magic squares $\mathbf{M}_{m}$ and $\mathbf{M}_{n},$ two order- $mn$
matrices can be formed by compounding as follows \cite{EGGE,ROGE}:%
\begin{equation}
\mathbf{A}_{mn}=\mathbf{E}_{m}\mathbf{\otimes M}_{n},\quad \mathbf{B}_{mn}=%
\mathbf{M}_{m}\mathbf{\otimes E}_{n},  \label{AB}
\end{equation}%
where the symbol $"\mathbf{\otimes "}$ denotes the tensor (Kronecker)
product \cite{HORN,MEYE}. On noting that%
\begin{align}
\mathbf{E}_{m}\mathbf{\otimes E}_{n}& =\mathbf{E}_{n}\mathbf{\otimes E}_{m}=%
\mathbf{E}_{mn},\quad \mathbf{E}_{n}^{2}=n\mathbf{E}_{n},\quad \limfunc{tr}%
\left[ \mathbf{E}_{n}\right] =n,  \notag \\
\mathbf{R}_{m}\mathbf{\otimes R}_{n}& =\mathbf{R}_{n}\mathbf{\otimes R}_{m}=%
\mathbf{R}_{mn},\quad \mathbf{R}_{n}^{2}=\mathbf{I}_{n},\quad \mathbf{R}_{n}%
\mathbf{E}_{n}=\mathbf{E}_{n}\mathbf{R}_{n}=\mathbf{E}_{n},  \label{EEn}
\end{align}%
and using formulas for the tensor product \cite{MEYE}, from (\ref{MnMagic}),
(\ref{AB}), and (\ref{EEn}), we find that%
\begin{align}
\mathbf{E}_{mn}\mathbf{A}_{mn}& =\left( \mathbf{E}_{m}\mathbf{\otimes E}%
_{n}\right) \left( \mathbf{E}_{m}\mathbf{\otimes M}_{n}\right) =\mathbf{E}%
_{m}^{2}\mathbf{\otimes }\left( \mathbf{E}_{n}\mathbf{M}_{n}\right) =m%
\mathbf{E}_{m}\mathbf{\otimes }\left( \mu _{n}\mathbf{E}_{n}\right) =m\mu
_{n}\mathbf{E}_{mn},  \notag \\
\limfunc{tr}\left[ \mathbf{A}_{mn}\right] & =\limfunc{tr}\left[ \mathbf{E}%
_{m}\mathbf{\otimes M}_{n}\right] =\limfunc{tr}\left[ \mathbf{E}_{m}\right] 
\limfunc{tr}\left[ \mathbf{M}_{n}\right] =m\mu _{n},  \label{ABmagic} \\
\limfunc{tr}\left[ \mathbf{R}_{mn}\mathbf{A}_{mn}\right] & =\limfunc{tr}%
\left[ \left( \mathbf{R}_{m}\mathbf{\otimes R}_{n}\right) \left( \mathbf{E}%
_{m}\mathbf{\otimes M}_{n}\right) \right] =\limfunc{tr}\left[ \left( \mathbf{%
R}_{m}\mathbf{E}_{m}\right) \mathbf{\otimes }\left( \mathbf{R}_{n}\mathbf{M}%
_{n}\right) \right] =m\mu _{n},  \notag
\end{align}%
and similarly for $\mathbf{B}_{mn}.$ Thus, $\mathbf{A}_{mn}$ and $\mathbf{B}%
_{mn}$ are unnatural magic squares with summation indices $m\mu _{n}$ and $%
n\mu _{m}$, respectively. Furthermore, we find that $\mathbf{A}_{mn}$ and $%
\mathbf{B}_{mn}$ commute since%
\begin{align}
\mathbf{A}_{mn}\mathbf{B}_{mn}& =\left( \mathbf{E}_{m}\mathbf{\otimes M}%
_{n}\right) \left( \mathbf{M}_{m}\mathbf{\otimes E}_{n}\right) =\left( 
\mathbf{E}_{m}\mathbf{M}_{m}\right) \mathbf{\otimes }\left( \mathbf{M}_{n}%
\mathbf{E}_{n}\right) =\mu _{m}\mu _{n}\mathbf{E}_{mn},  \notag \\
\mathbf{B}_{mn}\mathbf{A}_{mn}& =\left( \mathbf{M}_{m}\mathbf{\otimes E}%
_{n}\right) \left( \mathbf{E}_{m}\mathbf{\otimes M}_{n}\right) =\left( 
\mathbf{M}_{m}\mathbf{E}_{m}\right) \mathbf{\otimes }\left( \mathbf{E}_{n}%
\mathbf{M}_{n}\right) =\mu _{m}\mu _{n}\mathbf{E}_{mn}.  \label{ABcomm}
\end{align}%
In addition, it can be shown that each of the eight phases \cite{LOLY,NORD0}
of $\mathbf{A}_{mn}$ commutes with the corresponding phase of $\mathbf{B}%
_{mn}.$

We note that $\mathbf{A}_{mn}$ and $\mathbf{B}_{mn}$ are an orthogonal pair
since the $n$ by $n$ subsquares of $\mathbf{A}_{mn}$ contain $m^{2}$
replicas of $\mathbf{M}_{n},$ whereas the elements of the $n$ by $n$
subsquares of $\mathbf{B}_{mn}$ all are the same number as that of the
corresponding element of $\mathbf{M}_{m}$ (as seen in the examples below).
Thus, all combinations of two numbers from $\mathbf{A}_{mn}$ and $\mathbf{B}%
_{mn}$ occur once and only once, i.e. $\mathbf{A}_{mn}$ and $\mathbf{B}_{mn}$
are orthogonal. Therefore, we can form two natural magic squares by the
Euler composition formula \cite{PASL}\footnote{%
Pasles \cite{PASL} suggests that Benjamin Franklin used this formula prior
to Euler.} as follows:%
\begin{equation}
\mathbf{M}_{mn}^{A}=\mathbf{A}_{mn}+n^{2}\mathbf{B}_{mn},\qquad \mathbf{M}%
_{mn}^{B}=\mathbf{B}_{mn}+m^{2}\mathbf{A}_{mn}  \label{Mn2s}
\end{equation}%
which can be verified to satisfy the magic square conditions (\ref{MnMagic}%
), e.g.%
\begin{align}
\mathbf{E}_{mn}\mathbf{M}_{mn}^{A}& =\mathbf{E}_{mn}\mathbf{A}_{mn}+n^{2}%
\mathbf{E}_{mn}\mathbf{B}_{mn}=\left( m\mu _{n}+n^{3}\mu _{m}\right) \mathbf{%
E}_{mn}  \notag \\
& =\frac{1}{2}mn\left( m^{2}n^{2}-1\right) \mathbf{E}_{mn}=\mu _{mn}\mathbf{E%
}_{mn},  \label{EMn2} \\
\limfunc{tr}\mathbf{M}_{mn}^{A}& =\limfunc{tr}\mathbf{A}_{mn}+n^{2}\limfunc{%
tr}\mathbf{B}_{mn}=m\mu _{n}+n^{3}\mu _{m}=\mu _{mn}.  \notag
\end{align}%
We note that $\mathbf{M}_{mn}^{\left( A\right) }$ and $\mathbf{M}%
_{mn}^{\left( B\right) }$ commute since, by (\ref{Mn2s}) and (\ref{ABcomm})%
\begin{align}
\mathbf{M}_{mn}^{A}\mathbf{M}_{mn}^{B}& =\left( \mathbf{A}_{mn}+n^{2}\mathbf{%
B}_{mn}\right) \left( \mathbf{B}_{mn}+m^{2}\mathbf{A}_{mn}\right)   \notag \\
& =\mathbf{A}_{mn}\mathbf{B}_{mn}+m^{2}\mathbf{A}_{mn}\mathbf{A}_{mn}+n^{2}%
\mathbf{B}_{mn}\mathbf{B}_{mn}+m^{2}n^{2}\mathbf{B}_{mn}\mathbf{A}_{mn}, 
\notag \\
\mathbf{M}_{mn}^{B}\mathbf{M}_{mn}^{A}& =\left( \mathbf{B}_{mn}+m^{2}\mathbf{%
A}_{mn}\right) \left( \mathbf{A}_{mn}+n^{2}\mathbf{B}_{mn}\right)  \\
& =\mathbf{B}_{mn}\mathbf{A}_{mn}+n^{2}\mathbf{B}_{mn}\mathbf{B}_{mn}+m^{2}%
\mathbf{A}_{mn}\mathbf{A}_{mn}+m^{2}n^{2}\mathbf{A}_{mn}\mathbf{B}_{mn}, 
\notag \\
\mathbf{M}_{mn}^{A}\mathbf{M}_{mn}^{B}& =\mathbf{M}_{mn}^{B}\mathbf{M}%
_{mn}^{A}.  \notag
\end{align}

The foregoing compounding construction can be repeated using $\mathbf{M}%
_{mn}^{A}$ or $\mathbf{M}_{mn}^{B}$ in (\ref{AB}) to produce higher order
magic squares which again commute. Gigantic commuting magic squares can be
produced by repeated compounding as done by Chan and Loly \cite{CHAN}%
.\smallskip

\noindent \textbf{Regular Magic Squares. }In a \emph{regular} (associative)
matrix any two elements that are symmetric about the center element add to
the same number and in an odd-order regular matrix the center element is
one-half this number. The regularity condition on $\mathbf{M}_{n}$ can be
expressed as%
\begin{equation}
\mathbf{M}_{n}+\mathbf{R}_{n}\mathbf{M}_{n}\mathbf{R}_{n}=\frac{2\mu _{n}}{n}%
\mathbf{E}_{n}=\left( n^{2}-1\right) \mathbf{E}_{n},  \label{Rn}
\end{equation}%
where the factor $2\mu _{n}/n$ can be verified by taking the trace of this
equation. We wish to show that if $\mathbf{M}_{m}$ and $\mathbf{M}_{n}$ are
regular, then so are $\mathbf{A}_{mn}$, $\mathbf{B}_{mn}$, $\mathbf{M}%
_{mn}^{A},$ and $\mathbf{M}_{mn}^{B}$ as noted by Eggermont \cite{EGGE} for $%
\mathbf{M}_{mn}^{A}$. From (\ref{EEn}), (\ref{AB}), and (\ref{Rn}), we find
that%
\begin{align}
\mathbf{R}_{mn}\mathbf{A}_{mn}\mathbf{R}_{mn}+\mathbf{A}_{mn}& =\left( 
\mathbf{R}_{m}\mathbf{\otimes R}_{n}\right) \left( \mathbf{E}_{m}\mathbf{%
\otimes M}_{n}\right) \left( \mathbf{R}_{m}\mathbf{\otimes R}_{n}\right) +%
\mathbf{A}_{mn}  \notag \\
& =\left( \mathbf{R}_{m}\mathbf{E}_{m}\mathbf{R}_{m}\right) \mathbf{\otimes }%
\left( \mathbf{R}_{n}\mathbf{M}_{n}\mathbf{R}_{n}\right) +\mathbf{A}_{mn} 
\notag \\
& =\mathbf{E}_{m}\mathbf{\otimes }\left( \left( n^{2}-1\right) \mathbf{E}%
_{n}-\mathbf{M}_{n}\right) +\mathbf{E}_{m}\mathbf{\otimes M}_{n}
\label{AReg} \\
& =\left( n^{2}-1\right) \mathbf{E}_{mn}  \notag
\end{align}%
and similarly%
\begin{equation}
\mathbf{R}_{mn}\mathbf{B}_{mn}\mathbf{R}_{mn}+\mathbf{B}_{mn}=\left(
m^{2}-1\right) \mathbf{E}_{mn}  \label{BReg}
\end{equation}%
which are the regularity conditions for $\mathbf{A}_{mn}$ and $\mathbf{B}%
_{mn}$. From their definitions (\ref{Mn2s}), $\mathbf{M}_{mn}^{A}$ and $%
\mathbf{M}_{mn}^{B}$ are regular when $\mathbf{A}_{mn}$ and $\mathbf{B}_{mn}$
are regular as seen from 
\begin{align}
\mathbf{M}_{mn}^{A}+\mathbf{R}_{mn}\mathbf{M}_{mn}^{A}\mathbf{R}_{mn}&
=\left( \mathbf{A}_{mn}+n^{2}\mathbf{B}_{mn}\right) +\mathbf{R}_{mn}\left( 
\mathbf{A}_{mn}+n^{2}\mathbf{B}_{mn}\right) \mathbf{R}_{mn}  \notag \\
& =\mathbf{A}_{mn}+n^{2}\mathbf{B}_{mn}+\left( n^{2}-1\right) \mathbf{E}%
_{mn}-\mathbf{A}_{mn}+n^{2}\left( m^{2}-1\right) \mathbf{E}_{mn}-n^{2}%
\mathbf{B}_{mn}  \notag \\
& =\left( m^{2}n^{2}-1\right) \mathbf{E}_{mn}
\end{align}%
and similarly for $\mathbf{M}_{mn}^{B}$.

In addition, for a regular magic square $\mathbf{M}_{n}^{\left( 1\right) }$
and its $180^{\circ }$ rotation $\mathbf{M}_{n}^{\left( 2\right) },$ given by%
\begin{equation}
\mathbf{M}_{n}^{\left( 2\right) }=\mathbf{R}_{n}\mathbf{M}_{n}^{\left(
1\right) }\mathbf{R}_{n},  \label{M2nR}
\end{equation}%
it follows from (\ref{Rn}) that $\mathbf{M}_{n}^{\left( 1\right) }$ and $%
\mathbf{M}_{n}^{\left( 2\right) }$ commute and similarly for $\mathbf{M}%
_{m}^{\left( 1\right) }$ and $\mathbf{M}_{m}^{\left( 2\right) }$. These two
commuting duos can be used in equations of the form (\ref{AB}) and (\ref%
{Mn2s}) to form a quartet of mutually commuting regular magic squares.
Repeated compounding of these squares leads in an immense number of
commuting regular magic squares of increasing order. A class of pandiagonal
squares given in \cite{NORD0} also can be used to form commuting
duos.\smallskip 

\noindent \textbf{Pandiagonal Magic Squares.} In a \emph{pandiagonal} magic
square of order $n,$ all $2n$ diagonals, including broken ones in both
directions, sum to $\mu _{n}.$ It is known that a regular magic square $%
\mathbf{M}_{Rn}$ of doubly-even order $\left( n=4k,\ k=1,2,\ldots \right) $
can be transformed to a pandiagonal magic square $\mathbf{M}_{Pn}$ by the
Planck transformation \cite{PLAN,NORD0}. Thus, if $\mathbf{M}_{m}$ and $%
\mathbf{M}_{n}$ are regular and $mn$ is doubly-even, then $\mathbf{M}%
_{mn}^{A}$ and $\mathbf{M}_{mn}^{B}$ are regular (as shown above) and
doubly-even order. Therefore, they can be transformed to pandiagonal magic
squares which can be shown to commute.

It also is possible to compound commuting pandiagonal magic squares $\mathbf{%
M}_{mn}^{A}$ and $\mathbf{M}_{mn}^{B}$ directly from (\ref{AB}) and (\ref%
{Mn2s}) starting with pandiagonal magic squares $\mathbf{M}_{m}$ and $%
\mathbf{M}_{n}$ as noted by Eggermont \cite{EGGE} and carried out by Chan
and Loly \cite{CHAN} for $\mathbf{M}_{16}^{A}$. The pandiagonality of the
compounded matrix is established by them and is verified in an example below
for $\mathbf{M}_{16}^{A}$ and $\mathbf{M}_{16}^{B}$. Magic squares that are
both regular and pandiagonal are called \emph{ultra-magic} squares, with $5$
being the lowest order for their existence, leading to order-$25$ commuting
ultra-magic squares by the compounding construction.\smallskip 

\noindent \textbf{Special Case - Permutation}. Simplification is possible by
taking $\mathbf{M}_{m}=\mathbf{M}_{n}$ $\left( m=n\right) $ in (\ref{AB}).
In this case, as noted by Rogers, et.\negthinspace\ al \cite{ROGE},
interchange (shuffling) of rows and columns of $\mathbf{A}_{nn}$ leads to $%
\mathbf{B}_{nn}$ and vice versa. We find that this interchange can be
expressed as%
\begin{equation}
\mathbf{B}_{nn}=\mathbf{P}_{nn}\mathbf{A}_{nn}\mathbf{P}_{nn},\quad \mathbf{A%
}_{nn}=\mathbf{P}_{nn}\mathbf{B}_{nn}\mathbf{P}_{nn},  \label{Perm}
\end{equation}%
where $\mathbf{P}_{nn}$ is a symmetric permutation matrix that can be
written in block form as 
\begin{equation}
\mathbf{P}_{nn}=\left[ 
\begin{array}{cccc}
\mathbf{p}_{11} & \mathbf{p}_{12} & \cdots & \mathbf{p}_{1n} \\ 
\mathbf{p}_{21} & \mathbf{p}_{22} & \cdots & \mathbf{p}_{2n} \\ 
\vdots & \vdots & \ddots & \vdots \\ 
\mathbf{p}_{n1} & \mathbf{p}_{n2} & \cdots & \mathbf{p}_{nn}%
\end{array}%
\right]  \label{Pn2}
\end{equation}%
in which $\mathbf{p}_{ij}$ are order-$n$ matrices with element $\mathbf{p}%
_{ij}\left( j,i\right) =1$ and all other elements zero. It can be shown that 
\begin{equation}
\mathbf{R}_{nn}\mathbf{P}_{nn}\mathbf{R}_{nn}=\mathbf{P}_{nn},\quad \mathbf{P%
}_{nn}^{T}=\mathbf{P}_{nn}^{-1}=\mathbf{P}_{nn}.  \label{PRP}
\end{equation}%
Thus,(\ref{Perm}) is a \emph{similarity transformation }\cite{MEYE} and $%
\mathbf{A}_{nn}$ and $\mathbf{B}_{nn}$ are \emph{similar }matrices\emph{. }%
The formulas (\ref{Perm}), (\ref{Pn2}) and (\ref{PRP}) will be verified in
the examples below. From (\ref{Perm}) and (\ref{Mn2s}), $\mathbf{M}_{nn}^{A}$
and $\mathbf{M}_{nn}^{B}$ also are similar and are related by the
permutations%
\begin{equation}
\mathbf{M}_{nn}^{B}=\mathbf{P}_{nn}\mathbf{M}_{nn}^{A}\mathbf{P}_{nn},\qquad 
\mathbf{M}_{nn}^{A}=\mathbf{P}_{nn}\mathbf{M}_{nn}^{B}\mathbf{P}_{nn}.
\label{PermM}
\end{equation}

\noindent \textbf{Extension}. A huge number of pairs of commuting magic
squares can be constructed using various $\mathbf{M}_{n}$'s as the $m^{2}$
subsquares of \emph{generalized} $\mathbf{\tilde{A}}_{mn}$ and any $\mathbf{M%
}_{m}$ as the basis for $\mathbf{B}_{mn}$ in (\ref{AB}). Again, $\mathbf{%
\tilde{A}}_{mn}$ and $\mathbf{B}_{mn}$ are an orthogonal pair and they
commute since it can be shown that they satisfy (\ref{ABcomm}). Thus,
commuting $\mathbf{\tilde{M}}_{mn}^{A}$ and $\mathbf{\tilde{M}}_{mn}^{B}$
can be formed from them using (\ref{Mn2s}). However, since such an $\mathbf{%
\tilde{A}}_{mn}$ is not of the form (\ref{AB}), the general formulas for its
Jordan form and SVD and those of $\mathbf{\tilde{M}}_{mn}^{A}$ and $\mathbf{%
\tilde{M}}_{mn}^{B}$ (to be found next) do not apply.

\section{ Jordan Form}

We derive formulas for the Jordan-form matrices of $\mathbf{A}_{mn}$, $%
\mathbf{B}_{mn}$, $\mathbf{M}_{mn}^{A}$, and $\mathbf{M}_{mn}^{B}$ from
those of $\mathbf{E}_{m},\mathbf{E}_{n},\mathbf{M}_{m}$, and $\mathbf{M}_{n}$%
. To review from \cite{HORN,MEYE}, the Jordan form of a square matrix $%
\mathbf{M}$ is given by%
\begin{equation}
\mathbf{M=SJS}^{-1},\qquad \mathbf{MS=SJ,}  \label{LS}
\end{equation}%
where $\mathbf{S}$ is a matrix whose columns are the simple or generalized
eigenvectors $\mathbf{s}_{i}$ of $\mathbf{M}$ and $\mathbf{J}$ is the matrix
with zero elements except for eigenvalues $\lambda _{i}$ on the main
diagonal and ones on the diagonal above it corresponding to generalized
eigenvectors. For a generalized eigenvector $\mathbf{s}_{k}^{\left( i\right)
}$ with eigenvalue $\lambda _{i}$ of algebraic multiplicity $k,$ (\ref{LS})
leads to%
\begin{equation}
\left( \mathbf{M\,-\,}\lambda _{i}\mathbf{I}\right) ^{k}\mathbf{s}%
_{k}^{\left( i\right) }=\left( \mathbf{M\,-\,}\lambda _{i}\mathbf{I}\right)
^{k-1}\mathbf{s}_{k-1}^{\left( i\right) }=\mathbf{\ldots }=\left( \mathbf{%
M\,-\,}\lambda _{i}\mathbf{I}\right) \mathbf{s}_{1}^{\left( i\right) }=%
\mathbf{0.}  \label{Mlam}
\end{equation}%
For a simple eigenvector $\mathbf{s}_{i}$ with eigenvalue $\lambda _{i}$ of
algebraic multiplicity $1,$ (\ref{LS}) gives 
\begin{equation}
\left( \mathbf{M\,-\,}\lambda _{i}\mathbf{I}\right) \mathbf{s}_{i}=\mathbf{0.%
}
\end{equation}%
If all the eigenvectors of $\mathbf{M\,}$ are simple, then $\mathbf{J\equiv D%
}$ is diagonal.

The Jordan form of magic square matrices is studied extensively in \cite%
{LOLY,MATT,NORD0,NORD2,ROGE}. When $\mathbf{M}_{n}$ is a magic square
matrix, according to (\ref{MnMagic}), it has an all-ones eigenvector $%
\mathbf{s}_{1}$ with eigenvalue $\lambda _{1}=\mu _{n}.$ On applying $%
\mathbf{E}_{n}$ to (\ref{Mlam}), we find that%
\begin{equation}
\left( \mu _{n}\,\mathbf{-\,}\lambda _{i}\right) ^{k}\mathbf{E}_{n}\mathbf{s}%
_{k}^{\left( i\right) }=\left( \mu _{n}\,\mathbf{-\,}\lambda _{i}\right)
^{k-1}\mathbf{E}_{n}\mathbf{s}_{k-1}^{\left( i\right) }=\mathbf{\ldots }%
=\left( \mu _{n}\,\mathbf{-\,}\lambda _{i}\right) \mathbf{E}_{n}\mathbf{s}%
_{1}^{\left( i\right) }=\mathbf{0.}
\end{equation}%
Since it is known \cite{MATT} that $\left\vert \lambda _{i}\right\vert <\mu
_{n}$ $\left( _{i}\geq 2\right) $, it follows that 
\begin{equation}
\mathbf{E}_{n}\mathbf{s}_{1}=n\mu _{n}\mathbf{s}_{1},\quad \mathbf{E}_{n}%
\mathbf{s}_{i}=\mathbf{0},\quad i=2,3,\ldots .
\end{equation}%
The eigenvalues of $\mathbf{E}_{n}$ are $n,0,0,\ldots ,0$ as shown in \cite%
{HORN}. Thus, from (\ref{LS}) for $\mathbf{E}_{n},$ using the $\mathbf{S}%
_{n} $ matrix of $\mathbf{M}_{n},$ we have%
\begin{align}
\mathbf{E}_{n}\mathbf{S}_{n}\mathbf{\,}& \mathbf{=S}_{n}\mathbf{D}_{En} 
\notag \\
\mathbf{E}_{n}\left[ \mathbf{s}_{1},\mathbf{s}_{2},\ldots ,\mathbf{s}_{n}%
\right] & =\left[ \mathbf{s}_{1},\mathbf{s}_{2},\ldots ,\mathbf{s}_{n}\right]
\limfunc{diag}\left[ n,0,0,\ldots ,0\right] \\
\left[ n\mathbf{s}_{1},0,\ldots ,0\right] & =\left[ n\mathbf{s}_{1},0,\ldots
,0\right]  \notag
\end{align}%
which is correct, whence $\mathbf{S}_{n}$ is the eigenvector matrix for both 
$\mathbf{M}_{n}$ and $\mathbf{E}_{n}\mathbf{\ }$and their Jordan forms read%
\begin{equation}
\mathbf{M}_{n}=\mathbf{S}_{n}\mathbf{J}_{Mn}\mathbf{S}_{n}^{-1},\quad 
\mathbf{E}_{n}=\mathbf{S}_{n}\mathbf{D}_{En}\mathbf{S}_{n}^{-1}.  \label{EMS}
\end{equation}%
Since $\mathbf{E}_{n}$ is symmetric, it also has an orthogonal eigenvector
matrix which is not used here. When all the eigenvectors of $\mathbf{M}_{n}$
are simple, its eigenvalue matrix can be written as%
\begin{equation}
\mathbf{J}_{Mn}=\limfunc{diag}\left( \mu _{n},\lambda _{n2},\lambda
_{n3},\ldots ,\lambda _{nn}\right)  \label{DeDm}
\end{equation}%
and for generalized eigenvectors of $\mathbf{M}_{n}$ there are ones on the
diagonal above the main diagonal for their corresponding eigenvalues in $%
\mathbf{J}_{Mn}.$ Equations of the same form as the foregoing ones apply to $%
\mathbf{M}_{m}$ and $\mathbf{E}_{m}.$

Using a compounding technique given by Nordgren \cite{NORD2}, from (\ref{EMS}%
) and (\ref{AB}), we find that%
\begin{align}
\mathbf{A}_{mn}& =\mathbf{E}_{m}\mathbf{\otimes M}_{n}=\left( \mathbf{S}_{m}%
\mathbf{D}_{Em}\mathbf{S}_{m}^{-1}\right) \mathbf{\otimes }\left( \mathbf{S}%
_{n}\mathbf{J}_{Mn}\mathbf{S}_{n}^{-1}\right)   \notag \\
& =\left( \mathbf{S}_{m}\mathbf{\otimes S}_{n}\right) \left( \mathbf{D}_{Em}%
\mathbf{\otimes J}_{Mn}\right) \left( \mathbf{S}_{m}\mathbf{\otimes S}%
_{n}\right) ^{-1}=\mathbf{S}_{mn}\mathbf{J}_{mn}^{A}\mathbf{S}_{mn}^{-1}, 
\notag \\
\mathbf{B}_{mn}& =\mathbf{M}_{m}\mathbf{\otimes E}_{n}=\left( \mathbf{S}_{m}%
\mathbf{J}_{Mm}\mathbf{S}_{m}^{-1}\right) \mathbf{\otimes }\left( \mathbf{S}%
_{n}\mathbf{D}_{En}\mathbf{S}_{n}^{-1}\right)   \label{ABJord} \\
& =\left( \mathbf{S}_{nm}\mathbf{\otimes S}_{n}\right) \left( \mathbf{J}_{Mm}%
\mathbf{\otimes D}_{En}\right) \left( \mathbf{S}_{m}\mathbf{\otimes S}%
_{n}\right) ^{-1}=\mathbf{S}_{mn}\mathbf{J}_{mn}^{B}\mathbf{S}_{mn}^{-1}, 
\notag
\end{align}%
where\footnote{%
Rogers, et. \negthinspace al \cite{ROGE} derive similar formulas for
eigenvalues and eigenvectors in a somewhat different manner.}%
\begin{equation}
\mathbf{S}_{mn}=\mathbf{S}_{m}\mathbf{\otimes S}_{n},\quad \mathbf{J}%
_{mn}^{A}=\mathbf{D}_{Em}\mathbf{\otimes J}_{Mn},\quad \mathbf{J}_{mn}^{B}=%
\mathbf{J}_{Mm}\mathbf{\otimes D}_{En}.  \label{SmnD}
\end{equation}%
Then, with (\ref{DeDm}), it follows that the nonzero eigenvalues of $\mathbf{%
A}_{mn}$ and $\mathbf{B}_{mn}$ are%
\begin{align}
\mathbf{J}_{mn}^{A}& :n\mu _{m},n\lambda _{m2},n\lambda _{m3},\ldots
,n\lambda _{mn},  \notag \\
\mathbf{J}_{mn}^{B}& :m\mu _{n},m\lambda _{n2},m\lambda _{n3},\ldots
,m\lambda _{mn}.  \label{DADB}
\end{align}%
When $\mathbf{M}_{m}$ and/or $\mathbf{M}_{n}$ have generalized eigenvectors, 
$\mathbf{J}_{mn}^{A}$ and $\mathbf{J}_{mn}^{B}$ from (\ref{SmnD}) are not in
standard form but they can be brought there by modifying $\mathbf{S}_{mn}$
as indicated in the example below for $mn=12$. Furthermore, (\ref{Mn2s})
with (\ref{ABJord}) gives the Jordan form of $\mathbf{M}_{mn}^{A}$ and $%
\mathbf{M}_{mn}^{B}$ as%
\begin{equation}
\mathbf{M}_{mn}^{A}=\mathbf{S}_{mn}\mathbf{\mathbf{J}}_{Mmn}^{A}\mathbf{S}%
_{mn}^{-1},\quad \mathbf{M}_{mn}^{B}=\mathbf{S}_{mn}\mathbf{\mathbf{J}}%
_{Mmn}^{B}\mathbf{S}_{mn}^{-1},  \label{MJord}
\end{equation}%
where%
\begin{equation}
\mathbf{\mathbf{J}}_{Mmn}^{A}=\mathbf{J}_{mn}^{A}+n^{2}\mathbf{J}%
_{mn}^{B},\quad \mathbf{\mathbf{J}}_{Mmn}^{B}=\mathbf{J}_{mn}^{B}+m^{2}%
\mathbf{J}_{mn}^{A},  \label{Dn2}
\end{equation}%
and their eigenvalues can be expressed using (\ref{DADB}).\smallskip 

\noindent \textbf{Special Case}. In the special case where $\mathbf{M}_{m}=%
\mathbf{M}_{n}\ \left( m=n\right) ,$ according to (\ref{SmnD}), $\mathbf{A}%
_{nn}$ and $\mathbf{B}_{nn}$ have the same nonzero eigenvalues from (\ref%
{DADB}), namely%
\begin{equation}
n\mu _{n},n\lambda _{2},n\lambda _{3},\ldots ,n\lambda _{n},
\end{equation}%
but they appear in a different order in $\mathbf{J}_{nn}^{A}$ and $\mathbf{J}%
_{nn}^{B}.$ To see this, by (\ref{Perm}), (\ref{SmnD}), and (\ref{ABJord}),
we form%
\begin{equation}
\mathbf{A}_{nn}=\left( \mathbf{P}_{nn}\mathbf{S}_{nn}\mathbf{P}_{nn}\right)
\left( \mathbf{P}_{nn}\mathbf{J}_{nn}^{B}\mathbf{P}_{nn}\right) \left( 
\mathbf{P}_{nn}\mathbf{S}_{nn}\mathbf{P}_{nn}\right) ^{-1}=\mathbf{S}_{nn}%
\mathbf{J}_{nn}^{A}\mathbf{S}_{nn}^{-1},  \label{ABPn2}
\end{equation}%
and similarly for $\mathbf{B}_{nn},$ whence%
\begin{equation}
\mathbf{P}_{nn}\mathbf{S}_{nn}\mathbf{P}_{nn}=\mathbf{S}_{nn},\quad \mathbf{J%
}_{nn}^{A}=\mathbf{P}_{nn}\mathbf{J}_{nn}^{B}\mathbf{P}_{nn},\quad \mathbf{J}%
_{nn}^{B}=\mathbf{P}_{nn}\mathbf{J}_{nn}^{A}\mathbf{P}_{nn}  \label{PSP}
\end{equation}%
which confirms that $\mathbf{J}_{nn}^{A}$ and $\mathbf{J}_{nn}^{B}$ contain
the same eigenvalues and indicates their reordering. Furthermore, $\mathbf{M}%
_{nn}^{A}$ and $\mathbf{M}_{nn}^{B}$ also have the same nonzero eigenvalues,
namely%
\begin{equation}
n\left( 1+n^{2}\right) \mu _{n},n\lambda _{2},n\lambda _{3},\ldots ,n\lambda
_{n},n^{3}\lambda _{2},n^{3}\lambda _{3},\ldots ,n^{3}\lambda _{n}
\label{evMn2}
\end{equation}%
and equations of the form (\ref{PSP}) apply to $\mathbf{\mathbf{J}}%
_{Mnn}^{A} $ and $\mathbf{\mathbf{J}}_{Mnn}^{B}$.\smallskip

\section{Singular Value Decomposition}

We derive formulas for the matrices in the singular value decomposition
(SVD) of $\mathbf{A}_{mn}$, $\mathbf{B}_{mn}$, $\mathbf{M}_{mn}^{A}$, and $%
\mathbf{M}_{mn}^{B}$ in terms of those of $\mathbf{E}_{m},$ $\mathbf{E}_{n}$%
, $\mathbf{M}_{m}$, and $\mathbf{M}_{n}$. To review \cite{HORN,MEYE}, the
SVD of any real square matrix $\mathbf{M}$ is expressed as%
\begin{equation}
\mathbf{M=U\Sigma V}^{T},  \label{Msvd}
\end{equation}%
where $\mathbf{U}$ and $\mathbf{V}$ are orthogonal matrices, and $\mathbf{%
\Sigma }$ is a diagonal matrix with non-negative real numbers (the singular
values) on the diagonal.\footnote{%
The SVD also applies to complex matrices and rectangular matrices which are
not considered here.} It follows from (\ref{Msvd}) that%
\begin{equation}
\mathbf{MM}^{T}\mathbf{=U\Sigma }^{2}\mathbf{U}^{T},\quad \mathbf{M}^{T}%
\mathbf{M=V\mathbf{\Sigma }^{2}V}^{T}  \label{MMt}
\end{equation}%
which are Jordan forms of the symmetric, positive semi-definite matrices $%
\mathbf{MM}^{T}$ and $\mathbf{M}^{T}\mathbf{M}.$ In particular, the Jordan
form of $\mathbf{MM}^{T}$ can be used to determine $\mathbf{U}$ and $\mathbf{%
\Sigma }$ after which $\mathbf{V}$ can be determined from (\ref{Msvd}). If $%
\mathbf{\Sigma }$ is nonsingular, then (\ref{Msvd}) leads to%
\begin{equation}
\mathbf{V}=\mathbf{M}^{T}\mathbf{U\Sigma }^{-1}.  \label{VA}
\end{equation}%
If $\mathbf{\Sigma }$ is singular, then an alternate approach given by Meyer 
\cite{MEYE} applies.

The SVD of magic square matrices is studied in \cite{LOLY,NORD2,ROGE}. When $%
\mathbf{M}_{n}$ is a magic square matrix, the $\mathbf{U}_{n}$ and $\mathbf{V%
}_{n}$ matrices for $\mathbf{M}_{n}$ also apply to $\mathbf{E}_{n}$ as we
show next. In view of (\ref{MnMagic}), $\mathbf{M}_{n}$ and $\mathbf{M}%
_{n}^{T}$ have an eigenvalue $\mu _{n}$ with eigenvector $\mathbf{\mathbf{s}}%
_{1}$ composed of constant elements $c$, therefore%
\begin{gather}
\mathbf{M}_{n}\mathbf{s}_{1}=\mathbf{M}_{n}^{T}\mathbf{s}_{1}=\mu _{n}%
\mathbf{s}_{1},  \notag \\
\mathbf{M}_{n}\mathbf{\mathbf{M}}_{n}^{T}\mathbf{s}_{1}=\mu _{n}^{2}\mathbf{s%
}_{1},  \label{Mu1} \\
\mathbf{s}_{1}=c\mathbf{e}_{n}\mathbf{,\quad E}_{n}\mathbf{s}_{1}=nc\mathbf{e%
}_{n}\mathbf{,}  \notag
\end{gather}%
where $\mathbf{e}_{n}$ is the order-$n$ (column) vector with all elements
one. From (\ref{MMt}) and (\ref{Mu1}), we see that $\mathbf{s}_{1}$ also is
an eigenvector in $\mathbf{U}_{n}$ for the singular value $\mu $. Since $%
\mathbf{U}_{n}$ is orthogonal, $\mathbf{u}_{1}=\mathbf{s}_{1}$ must be a
unit vector, hence%
\begin{equation}
\mathbf{u}_{1}^{T}\mathbf{u}_{1}=c^{2}\mathbf{e}_{n}^{T}\mathbf{e}%
_{n}=c^{2}n=1,\quad \therefore \mathbf{u}_{1}=\frac{\sqrt{n}}{n}\mathbf{e}%
_{n}\mathbf{,\quad E}_{n}\mathbf{u_{1}}=\sqrt{n}\mathbf{e}_{n}\mathbf{.}
\label{u1}
\end{equation}%
The remaining eigenvectors in $\mathbf{U}_{n}$ namely $\mathbf{u}_{2},%
\mathbf{u}_{3},\ldots ,\mathbf{u}_{n}$ for singular values $\sigma _{2}$,$%
\sigma _{3}$,$\ldots $,$\sigma _{n}$ in $\mathbf{\Sigma }_{Mn}\mathbf{,}$
according to (\ref{MMt}), must satisfy 
\begin{equation}
\left( \mathbf{\mathbf{M}}_{n}\mathbf{\mathbf{M}}_{n}^{T}-\sigma _{i}^{2}%
\mathbf{I}_{n}\right) \mathbf{u}_{i}=0,\quad i=2,3,\ldots ,n.  \label{MMtui}
\end{equation}%
Application of $\mathbf{E}_{n}$ to this equation results in 
\begin{equation}
\left( \mu ^{2}-\sigma _{i}^{2}\right) \mathbf{E}_{n}\mathbf{u}_{i}=0,\quad
i=2,3,\ldots ,n  \label{m2ui}
\end{equation}%
and, since it is known \cite{MATT} that $\sigma _{i}^{2}<\mu ^{2}$ $\left(
i\geq 2\right) $, it follows that%
\begin{equation}
\mathbf{E}_{n}\mathbf{u}_{i}=\mathbf{0},\quad i=2,3,\ldots ,n.  \label{Eui}
\end{equation}%
A similar argument holds for the singular vectors $\mathbf{v}_{i}$ of $%
\mathbf{V}_{n}$ and we may write (in block form)%
\begin{gather}
\mathbf{U}_{n}=\left[ 
\begin{array}{cccc}
\frac{\sqrt{n}}{n}\mathbf{e}_{n} & \mathbf{u}_{2} & \ldots  & \mathbf{u}_{n}%
\end{array}%
\right] ,\quad \mathbf{V}_{n}=\left[ 
\begin{array}{cccc}
\frac{\sqrt{n}}{n}\mathbf{e}_{n} & \mathbf{v}_{2} & \ldots  & \mathbf{v}_{n}%
\end{array}%
\right] ,  \notag \\
\mathbf{E}_{n}\mathbf{U}_{n}=\mathbf{E}_{n}\mathbf{V}_{n}=\left[ 
\begin{array}{cccc}
\sqrt{n}\mathbf{e}_{n} & \mathbf{0} & \ldots  & \mathbf{0}%
\end{array}%
\right] ,  \label{UVE}
\end{gather}%
where $\mathbf{v}_{2},\mathbf{v}_{3},\ldots \mathbf{v}_{n}$ remain to be
determined from (\ref{Msvd}) as already noted. From the SVD for $\mathbf{E}%
_{n}$ with $\mathbf{U}_{n}$ and $\mathbf{V}_{n}$ from (\ref{UVE}), we have 
\begin{equation}
\mathbf{\Sigma }_{En}=\mathbf{U}_{n}^{T}\mathbf{E}_{n}\mathbf{V}_{n}=\left[ 
\begin{array}{c}
\frac{\sqrt{n}}{n}\mathbf{e}_{n}^{T} \\ 
\mathbf{u}_{2}^{T} \\ 
\vdots  \\ 
\mathbf{u}_{n}^{T}%
\end{array}%
\right] \left[ 
\begin{array}{cccc}
\sqrt{n}\mathbf{e}_{n} & \mathbf{0} & \ldots  & \mathbf{0}%
\end{array}%
\right] =\limfunc{diag}\left[ n,0,0,\ldots ,0\right]   \label{UtEV}
\end{equation}%
which is correct. Therefore%
\begin{equation}
\mathbf{E}_{n}=\mathbf{U}_{n}\mathbf{\Sigma }_{En}\mathbf{V}_{n}^{T},\quad 
\mathbf{\Sigma }_{En}=\limfunc{diag}\left[ n,0,0,\ldots ,0\right] ,
\label{Esvd}
\end{equation}%
i.e. $\mathbf{E}_{n}$ has the same singular-value matrices $\mathbf{U}_{n}$
and $\mathbf{V}_{n}$ as $\mathbf{M}_{n}$.

Using a compounding technique given by Nordgren \cite{NORD2}, by (\ref{AB}),
(\ref{Msvd}), and (\ref{Esvd}), we have%
\begin{align}
\mathbf{A}_{mn}& =\mathbf{E}_{m}\mathbf{\otimes M}_{n}=\left( \mathbf{U}_{m}%
\mathbf{\Sigma }_{Em}\mathbf{V}_{m}^{T}\right) \mathbf{\otimes }\left( 
\mathbf{U}_{n}\mathbf{\Sigma }_{Mn}\mathbf{V}_{n}^{T}\right)   \notag \\
& =\left( \mathbf{U}_{m}\mathbf{\otimes U}_{n}\right) \left( \mathbf{\Sigma }%
_{Em}\mathbf{\otimes \Sigma }_{Mn}\right) \left( \mathbf{V}_{m}\mathbf{%
\otimes V}_{n}\right) ^{T}=\mathbf{U}_{mn}\mathbf{\Sigma }_{mn}^{A}\mathbf{V}%
_{mn}^{T},  \notag \\
\mathbf{B}_{mn}& =\mathbf{M}_{m}\mathbf{\otimes E}_{n}=\left( \mathbf{U}_{m}%
\mathbf{\Sigma }_{Mm}\mathbf{V}_{m}^{T}\right) \mathbf{\otimes }\left( 
\mathbf{U}_{n}\mathbf{\mathbf{\Sigma }}_{En}\mathbf{V}_{n}^{T}\right) 
\label{ABn2} \\
& =\left( \mathbf{U}_{m}\mathbf{\otimes U}_{n}\right) \left( \mathbf{\Sigma }%
_{Mm}\mathbf{\otimes \Sigma }_{En}\right) \left( \mathbf{V}_{m}\mathbf{%
\otimes V}_{n}\right) ^{T}=\mathbf{U}_{mn}\mathbf{\Sigma }_{mn}^{B}\mathbf{V}%
_{mn}^{T},  \notag
\end{align}%
where\footnote{%
Rogers, et. \negthinspace al \cite{ROGE} also give formulas for the SVD of
compound matrices.}%
\begin{align}
\mathbf{U}_{mn}& =\mathbf{U}_{m}\mathbf{\otimes U}_{n},\quad \mathbf{V}_{mn}=%
\mathbf{V}_{m}\mathbf{\otimes V}_{n},  \notag \\
\mathbf{\Sigma }_{mn}^{A}& =\mathbf{\Sigma }_{Em}\mathbf{\otimes \Sigma }%
_{Mn}=m\limfunc{diag}\left[ \mu _{n},\sigma _{n2},\sigma _{n3},\ldots
,\sigma _{nn},0,\ldots ,0\right] ,  \label{SVDAB} \\
\mathbf{\Sigma }_{mn}^{B}& =\mathbf{\Sigma }_{Mm}\mathbf{\otimes \Sigma }%
_{En}=n\limfunc{diag}\left[ \mu _{m},0,\ldots ,0,\sigma _{m2},0,\ldots
,0,\sigma _{m3},0,\ldots ,0,\ldots ,\sigma _{mm},0,\ldots 0\right] .  \notag
\end{align}%
Thus, $\mathbf{A}_{mn}$ and $\mathbf{B}_{mn}$ have the same $\mathbf{U}_{mn}$
and $\mathbf{V}_{mn}$ and their SVD's are given by (\ref{ABn2}).
Furthermore, it follows from (\ref{Dn2}) and (\ref{ABn2}) that the SVD's of $%
\mathbf{M}_{mn}^{A}$ and $\mathbf{M}_{mn}^{B}$ are%
\begin{equation}
\mathbf{M}_{mn}^{A}=\mathbf{U}_{mn}\mathbf{\Sigma }_{Mmn}^{A}\mathbf{V}%
_{mn}^{T},\quad \mathbf{M}_{mn}^{B}=\mathbf{U}_{mn}\mathbf{\Sigma }_{Mmn}^{B}%
\mathbf{V}_{mn}^{T},  \label{SVDMAB}
\end{equation}%
where%
\begin{equation}
\mathbf{\Sigma }_{Mmn}^{A}=\mathbf{\Sigma }_{mn}^{A}+m^{2}\mathbf{\Sigma }%
_{mn}^{B},\quad \mathbf{\Sigma }_{Mmn}^{B}=\mathbf{\Sigma }_{mn}^{B}+n^{2}%
\mathbf{\Sigma }_{mn}^{A}.  \label{SigAB}
\end{equation}%
\smallskip 

\noindent \textbf{Special Case}. In the special case where $\mathbf{M}_{m}=%
\mathbf{M}_{n}\ \left( m=n\right) ,$ according to (\ref{SVDAB}), $\mathbf{A}%
_{nn}$ and $\mathbf{B}_{nn}$ have the same singular values but they are in a
different order in $\mathbf{\Sigma }_{nn}^{A}$ and $\mathbf{\Sigma }%
_{nn}^{B}.$ To examine this, by (\ref{Perm}) and (\ref{ABn2}), we form 
\begin{align}
\mathbf{A}_{nn}& =\left( \mathbf{P}_{nn}\mathbf{U}_{nn}\mathbf{P}%
_{nn}\right) \left( \mathbf{P}_{nn}\mathbf{\Sigma }_{nn}^{B}\mathbf{P}%
_{nn}\right) \left( \mathbf{P}_{nn}\mathbf{V}_{nn}^{T}\mathbf{P}_{nn}\right)
^{-1}=\mathbf{U}_{nn}\mathbf{\Sigma }_{nn}^{A}\mathbf{V}_{nn}^{T},  \notag \\
\mathbf{B}_{nn}& =\left( \mathbf{P}_{nn}\mathbf{U}_{nn}\mathbf{P}%
_{nn}\right) \left( \mathbf{P}_{nn}\mathbf{\Sigma }_{nn}^{A}\mathbf{P}%
_{nn}\right) \left( \mathbf{P}_{nn}\mathbf{V}_{nn}^{T}\mathbf{P}_{nn}\right)
^{-1}=\mathbf{U}_{nn}\mathbf{\Sigma }_{nn}^{B}\mathbf{V}_{nn}^{T},
\label{PAP}
\end{align}%
whence%
\begin{gather}
\mathbf{P}_{nn}\mathbf{U}_{nn}\mathbf{P}_{nn}=\mathbf{U}_{nn},\quad \mathbf{P%
}_{nn}\mathbf{V}_{nn}\mathbf{P}_{nn}=\mathbf{V}_{nn},  \notag \\
\mathbf{\Sigma }_{nn}^{A}=\mathbf{P}_{nn}\mathbf{\Sigma }_{nn}^{B}\mathbf{P}%
_{nn},\quad \mathbf{\Sigma }_{nn}^{B}=\mathbf{P}_{nn}\mathbf{\Sigma }%
_{nn}^{A}\mathbf{P}_{nn}  \label{PUP}
\end{gather}%
which indicates the reordering of the same singular values in $\mathbf{%
\Sigma }_{nn}^{A}$ and $\mathbf{\Sigma }_{nn}^{B}$. The SVD's of $\mathbf{M}%
_{nn}^{A}$ and $\mathbf{M}_{nn}^{B}$ are still given by (\ref{SVDMAB}) with
singular value matrices from (\ref{SigAB}). By (\ref{SigAB}) and (\ref{PUP}%
), we have 
\begin{equation}
\mathbf{\Sigma }_{Mnn}^{A}=\mathbf{P}_{nn}\mathbf{\Sigma }_{Mnn}^{B}\mathbf{P%
}_{nn},\quad \mathbf{\Sigma }_{Mnn}^{B}=\mathbf{P}_{nn}\mathbf{\Sigma }%
_{Mnn}^{A}\mathbf{P}_{nn}  \label{Sign2P}
\end{equation}%
which confirms that $\mathbf{\Sigma }_{Mnn}^{A}$ and $\mathbf{\Sigma }%
_{Mnn}^{B}$ contain the same singular values and indicates their reordering.
Next, numerical examples are given that illustrate and confirm the foregoing
theoretical results.

\section{Examples}

\noindent \textbf{Order 9.} We construct two order-9, commuting, regular,
magic squares by compounding. Following Rogers, et.\negthinspace\ al \cite%
{ROGE}, we start with the order-$3$ Lo-Shu regular magic square $\mathbf{M}%
_{3}$ and the all-ones square $\mathbf{E}_{3},$ namely%
\begin{equation}
\mathbf{M}_{3}=\left[ 
\begin{array}{ccc}
3 & 8 & 1 \\ 
2 & 4 & 6 \\ 
7 & 0 & 5%
\end{array}%
\right] ,\quad \mathbf{E}_{3}=\left[ 
\begin{array}{ccc}
1 & 1 & 1 \\ 
1 & 1 & 1 \\ 
1 & 1 & 1%
\end{array}%
\right] .  \label{MLS}
\end{equation}%
We compound these matrices according to (\ref{AB}) to form%
\begin{equation}
\mathbf{A}_{9}=\left[ 
\begin{array}{ccccccccc}
3 & 8 & 1 & 3 & 8 & 1 & 3 & 8 & 1 \\ 
2 & 4 & 6 & 2 & 4 & 6 & 2 & 4 & 6 \\ 
7 & 0 & 5 & 7 & 0 & 5 & 7 & 0 & 5 \\ 
3 & 8 & 1 & 3 & 8 & 1 & 3 & 8 & 1 \\ 
2 & 4 & 6 & 2 & 4 & 6 & 2 & 4 & 6 \\ 
7 & 0 & 5 & 7 & 0 & 5 & 7 & 0 & 5 \\ 
3 & 8 & 1 & 3 & 8 & 1 & 3 & 8 & 1 \\ 
2 & 4 & 6 & 2 & 4 & 6 & 2 & 4 & 6 \\ 
7 & 0 & 5 & 7 & 0 & 5 & 7 & 0 & 5%
\end{array}%
\right] ,\ \ \mathbf{B}_{9}=\left[ 
\begin{array}{ccccccccc}
3 & 3 & 3 & 8 & 8 & 8 & 1 & 1 & 1 \\ 
3 & 3 & 3 & 8 & 8 & 8 & 1 & 1 & 1 \\ 
3 & 3 & 3 & 8 & 8 & 8 & 1 & 1 & 1 \\ 
2 & 2 & 2 & 4 & 4 & 4 & 6 & 6 & 6 \\ 
2 & 2 & 2 & 4 & 4 & 4 & 6 & 6 & 6 \\ 
2 & 2 & 2 & 4 & 4 & 4 & 6 & 6 & 6 \\ 
7 & 7 & 7 & 0 & 0 & 0 & 5 & 5 & 5 \\ 
7 & 7 & 7 & 0 & 0 & 0 & 5 & 5 & 5 \\ 
7 & 7 & 7 & 0 & 0 & 0 & 5 & 5 & 5%
\end{array}%
\right] ,
\end{equation}%
which are unnatural, regular, magic squares that commute. Since $\mathbf{A}%
_{9}$ and $\mathbf{B}_{9}$ are an orthogonal pair, two commuting, regular,
magic squares can be formed from $\mathbf{A}_{9}$ and $\mathbf{B}_{9}$
according to (\ref{Mn2s}) as%
\begin{gather}
\mathbf{M}_{9}^{A}=\hspace{2.3in}\mathbf{M}_{9}^{B}=  \notag \\
\left[ 
\begin{array}{ccccccccc}
30 & 35 & 28 & 75 & 80 & 73 & 12 & 17 & 10 \\ 
29 & 31 & 33 & 74 & 76 & 78 & 11 & 13 & 15 \\ 
34 & 27 & 32 & 79 & 72 & 77 & 16 & 9 & 14 \\ 
21 & 26 & 19 & 39 & 44 & 37 & 57 & 62 & 55 \\ 
20 & 22 & 24 & 38 & 40 & 42 & 56 & 58 & 60 \\ 
25 & 18 & 23 & 43 & 36 & 41 & 61 & 54 & 59 \\ 
66 & 71 & 64 & 3 & 8 & 1 & 48 & 53 & 46 \\ 
65 & 67 & 69 & 2 & 4 & 6 & 47 & 49 & 51 \\ 
70 & 63 & 68 & 7 & 0 & 5 & 52 & 45 & 50%
\end{array}%
\right] ,\quad \left[ 
\begin{array}{ccccccccc}
30 & 75 & 12 & 35 & 80 & 17 & 28 & 73 & 10 \\ 
21 & 39 & 57 & 26 & 44 & 62 & 19 & 37 & 55 \\ 
66 & 3 & 48 & 71 & 8 & 53 & 64 & 1 & 46 \\ 
29 & 74 & 11 & 31 & 76 & 13 & 33 & 78 & 15 \\ 
20 & 38 & 56 & 22 & 40 & 58 & 24 & 42 & 60 \\ 
65 & 2 & 47 & 67 & 4 & 49 & 69 & 6 & 51 \\ 
34 & 79 & 16 & 27 & 72 & 9 & 32 & 77 & 14 \\ 
25 & 43 & 61 & 18 & 36 & 54 & 23 & 41 & 59 \\ 
70 & 7 & 52 & 63 & 0 & 45 & 68 & 5 & 50%
\end{array}%
\right] .  \label{M911}
\end{gather}%
As noted in \cite{ROGE}, $\mathbf{M}_{9}^{A}$ was known prior to 1000 AD and 
$\mathbf{M}_{9}^{B}$ dates to 1275 AD.

In addition, as noted in Section 2, a huge number of pairs of commuting
magic squares $\mathbf{\tilde{M}}_{9}^{A}$ and $\mathbf{\tilde{M}}_{9}^{B}$
can be constructed according to (\ref{Mn2s}) using various combinations of
the eight phases of $\mathbf{M}_{3}$ as the nine subsquares of generalized $%
\mathbf{\tilde{A}}_{9}$ and any phase of $\mathbf{M}_{3}$ as the basis for $%
\mathbf{B}_{9}$ in (\ref{AB}). A regular $\mathbf{\tilde{A}}_{9}$ results
from using any five phases of $\mathbf{M}_{3}$ as the nine subsquares of $%
\mathbf{\tilde{A}}_{9}$ placed in a regular block pattern, e.g.%
\begin{equation}
\mathbf{\tilde{A}}_{9}=\left[ 
\begin{array}{ccccccccc}
5 & 6 & 1 & 3 & 2 & 7 & 3 & 8 & 1 \\ 
0 & 4 & 8 & 8 & 4 & 0 & 2 & 4 & 6 \\ 
7 & 2 & 3 & 1 & 6 & 5 & 7 & 0 & 5 \\ 
1 & 6 & 5 & 7 & 2 & 3 & 1 & 6 & 5 \\ 
8 & 4 & 0 & 0 & 4 & 8 & 8 & 4 & 0 \\ 
3 & 2 & 7 & 5 & 6 & 1 & 3 & 2 & 7 \\ 
3 & 8 & 1 & 3 & 2 & 7 & 5 & 6 & 1 \\ 
2 & 4 & 6 & 8 & 4 & 0 & 0 & 4 & 8 \\ 
7 & 0 & 5 & 1 & 6 & 5 & 7 & 2 & 3%
\end{array}%
\right] 
\end{equation}
It is easy to see that $\mathbf{\tilde{A}}_{9}$ and $\mathbf{B}_{9}$ are
orthogonal and they commute since%
\begin{equation}
\mathbf{\tilde{A}}_{9}\mathbf{B}_{9}=\mathbf{B}_{9}\mathbf{\tilde{A}}_{9}=144%
\mathbf{E}_{9}=\left( \mu _{9}\right) ^{2}\mathbf{E}_{9}.
\end{equation}%
Thus, as noted in Section 2, commuting regular $\mathbf{\tilde{M}}_{9}^{A}$
and $\mathbf{\tilde{M}}_{9}^{B}$ can be formed from them using (\ref{Mn2s}).

The permutation matrix $\mathbf{P}_{9}$ that connects $\mathbf{A}_{9}$ and $%
\mathbf{B}_{9}$ according to (\ref{Perm}) is given by (\ref{Pn2}) as 
\begin{equation}
\mathbf{P}_{9}=\left[ 
\begin{array}{ccccccccc}
1 & 0 & 0 & 0 & 0 & 0 & 0 & 0 & 0 \\ 
0 & 0 & 0 & 1 & 0 & 0 & 0 & 0 & 0 \\ 
0 & 0 & 0 & 0 & 0 & 0 & 1 & 0 & 0 \\ 
0 & 1 & 0 & 0 & 0 & 0 & 0 & 0 & 0 \\ 
0 & 0 & 0 & 0 & 1 & 0 & 0 & 0 & 0 \\ 
0 & 0 & 0 & 0 & 0 & 0 & 0 & 1 & 0 \\ 
0 & 0 & 1 & 0 & 0 & 0 & 0 & 0 & 0 \\ 
0 & 0 & 0 & 0 & 0 & 1 & 0 & 0 & 0 \\ 
0 & 0 & 0 & 0 & 0 & 0 & 0 & 0 & 1%
\end{array}%
\right] .  \label{P9}
\end{equation}%
$\allowbreak $It can be verified that $\mathbf{P}_{9}$ satisfies (\ref{PRP})
and connects $\mathbf{A}_{9}$ to $\mathbf{B}_{9}$ and $\mathbf{M}_{9}^{A}$
to $\mathbf{M}_{9}^{B}$ according to (\ref{Perm}) and (\ref{PermM}). Thus,
they are pairs of similar matrices as noted in Section 2.

In order to see why $\mathbf{P}_{9}$ works and can be generalized to higher
order $nn,$ we consider a general $\mathbf{M}_{3}$ compounded with $\mathbf{E%
}_{3}$ written in block form as%
\begin{equation}
\mathbf{B}_{9}=\mathbf{M}_{3}\mathbf{\otimes E}_{3}=\left[ 
\begin{array}{ccc}
\mathbf{b}_{11} & \mathbf{b}_{12} & \mathbf{b}_{13} \\ 
\mathbf{b}_{21} & \mathbf{b}_{22} & \mathbf{b}_{23} \\ 
\mathbf{b}_{31} & \mathbf{b}_{32} & \mathbf{b}_{33}%
\end{array}%
\right] ,\quad \mathbf{M}_{3}=\left[ 
\begin{array}{ccc}
m_{11} & m_{12} & m_{13} \\ 
m_{21} & m_{22} & m_{23} \\ 
m_{31} & m_{32} & m_{33}%
\end{array}%
\right] ,
\end{equation}%
where $\mathbf{b}_{ij}$ is a block that has all elements $m_{ij}.$ Then, the
permutation (\ref{Perm}) of $\mathbf{B}_{9}$ can be written as%
\begin{equation}
\mathbf{P}_{9}\mathbf{B}_{9}\mathbf{P}_{9}=\left[ 
\begin{array}{ccc}
\mathbf{p}_{11} & \mathbf{p}_{12} & \mathbf{p}_{13} \\ 
\mathbf{p}_{21} & \mathbf{p}_{22} & \mathbf{p}_{23} \\ 
\mathbf{p}_{31} & \mathbf{p}_{32} & \mathbf{p}_{33}%
\end{array}%
\right] \left[ 
\begin{array}{ccc}
\mathbf{b}_{11} & \mathbf{b}_{12} & \mathbf{b}_{13} \\ 
\mathbf{b}_{21} & \mathbf{b}_{22} & \mathbf{b}_{23} \\ 
\mathbf{b}_{31} & \mathbf{b}_{32} & \mathbf{b}_{33}%
\end{array}%
\right] \left[ 
\begin{array}{ccc}
\mathbf{p}_{11} & \mathbf{p}_{12} & \mathbf{p}_{13} \\ 
\mathbf{p}_{21} & \mathbf{p}_{22} & \mathbf{p}_{23} \\ 
\mathbf{p}_{31} & \mathbf{p}_{32} & \mathbf{p}_{33}%
\end{array}%
\right]
\end{equation}%
and on carrying out the matrix multiplication we find that%
\begin{equation}
\mathbf{P}_{9}\mathbf{B}_{9}\mathbf{P}_{9}=\left[ 
\begin{array}{ccc}
\mathbf{M}_{3} & \mathbf{M}_{3} & \mathbf{M}_{3} \\ 
\mathbf{M}_{3} & \mathbf{M}_{3} & \mathbf{M}_{3} \\ 
\mathbf{M}_{3} & \mathbf{M}_{3} & \mathbf{M}_{3}%
\end{array}%
\right] =\mathbf{E}_{3}\mathbf{\otimes M}_{3}=\mathbf{A}_{9}.  \label{P9B9}
\end{equation}%
As an example of this matrix multiplication, the element in the first row,
second column of $\mathbf{P}_{9}\mathbf{B}_{9}\mathbf{P}_{9}$ is given by%
\begin{align}
& \mathbf{p}_{11}\mathbf{b}_{11}\mathbf{p}_{12}+\mathbf{p}_{12}\mathbf{b}%
_{21}\mathbf{p}_{12}+\mathbf{p}_{13}\mathbf{b}_{31}\mathbf{p}_{12}+\mathbf{p}%
_{11}\mathbf{b}_{12}\mathbf{p}_{22}+\mathbf{p}_{12}\mathbf{b}_{22}\mathbf{p}%
_{22}  \notag \\
& +\mathbf{p}_{13}\mathbf{b}_{32}\mathbf{p}_{22}+\mathbf{p}_{11}\mathbf{b}%
_{13}\mathbf{p}_{32}+\mathbf{p}_{12}\mathbf{b}_{23}\mathbf{p}_{32}+\mathbf{p}%
_{13}\mathbf{b}_{33}\mathbf{p}_{32}  \notag \\
& =\left[ 
\begin{array}{ccc}
m_{11} & 0 & 0 \\ 
0 & 0 & 0 \\ 
0 & 0 & 0%
\end{array}%
\right] +\left[ 
\begin{array}{ccc}
0 & 0 & 0 \\ 
m_{21} & 0 & 0 \\ 
0 & 0 & 0%
\end{array}%
\right] +\left[ 
\begin{array}{ccc}
0 & 0 & 0 \\ 
0 & 0 & 0 \\ 
m_{31} & 0 & 0%
\end{array}%
\right] +\left[ 
\begin{array}{ccc}
0 & m_{12} & 0 \\ 
0 & 0 & 0 \\ 
0 & 0 & 0%
\end{array}%
\right]  \notag \\
& +\left[ 
\begin{array}{ccc}
0 & 0 & 0 \\ 
0 & m_{22} & 0 \\ 
0 & 0 & 0%
\end{array}%
\right] +\left[ 
\begin{array}{ccc}
0 & 0 & 0 \\ 
0 & 0 & 0 \\ 
0 & m_{32} & 0%
\end{array}%
\right] +\left[ 
\begin{array}{ccc}
0 & 0 & m_{13} \\ 
0 & 0 & 0 \\ 
0 & 0 & 0%
\end{array}%
\right] +\left[ 
\begin{array}{ccc}
0 & 0 & 0 \\ 
0 & 0 & m_{23} \\ 
0 & 0 & 0%
\end{array}%
\right] \\
& +\left[ 
\begin{array}{ccc}
0 & 0 & 0 \\ 
0 & 0 & 0 \\ 
0 & 0 & m_{33}%
\end{array}%
\right] =\left[ 
\begin{array}{ccc}
m_{11} & m_{12} & m_{13} \\ 
m_{21} & m_{22} & m_{23} \\ 
m_{31} & m_{32} & m_{33}%
\end{array}%
\right] =\mathbf{M}_{3}.  \notag
\end{align}%
In view of (\ref{P9B9}), we have verified (\ref{Perm}) for $n=3.$ It should
be clear that a similar verification of (\ref{Perm}) applies for higher
orders. Also, (\ref{PRP}) can be verified in a similar manner.

Next, we construct the Jordan-form matrices of $\mathbf{A}_{9},\ \mathbf{B}%
_{9},\ \mathbf{M}_{9}^{A},$ and $\mathbf{M}_{9}^{B}$ from the following
Jordan-form matrices of $\mathbf{M}_{3}$ and $\mathbf{E}_{3}$:%
\begin{align}
\mathbf{S}_{3}& =\left[ 
\begin{array}{ccc}
1 & 8+i\sqrt{6} & 8-i\sqrt{6} \\ 
1 & -4+2i\sqrt{6} & -4-2i\sqrt{6} \\ 
1 & -4-3i\sqrt{6} & -4+3i\sqrt{6}%
\end{array}%
\right] ,  \label{Jord3} \\
\mathbf{D}_{M3}& =\limfunc{diag}\left[ 12,2i\sqrt{6},-2i\sqrt{6}\right] , 
\notag \\
\mathbf{D}_{E3}& =\limfunc{diag}\left[ 3,0,0\right] ,  \notag
\end{align}%
where all the eigenvectors in $\mathbf{S}_{3}$ are simple. As noted in
Section 2, the eigenvector matrix $\mathbf{S}_{3}$ for $\mathbf{M}_{3}$ is
also an eigenvector matrix for $\mathbf{E}_{3}.$ Also, the eigenvalues of
the regular matrix $\mathbf{M}_{3}$ are $\mu _{3}$ and a $\pm $pair. By (\ref%
{Jord3}) and (\ref{SmnD}), we obtain the following Jordan-form matrices for $%
\mathbf{A}_{9}$ and $\mathbf{B}_{9}$:%
\begin{align}
\mathbf{S}_{9}& =\left[ 
\begin{array}{ccccc}
1 & 8+i\sqrt{6} & 8-i\sqrt{6} & 8+i\sqrt{6} & 58+16i\sqrt{6} \\ 
1 & -4+2i\sqrt{6} & -4-2i\sqrt{6} & 8+i\sqrt{6} & -44+12i\sqrt{6} \\ 
1 & -4-3i\sqrt{6} & -4+3i\sqrt{6} & 8+i\sqrt{6} & -14-28i\sqrt{6} \\ 
1 & 8+i\sqrt{6} & 8-i\sqrt{6} & -4+2i\sqrt{6} & -44+12i\sqrt{6} \\ 
1 & -4+2i\sqrt{6} & -4-2i\sqrt{6} & -4+2i\sqrt{6} & -8-16i\sqrt{6} \\ 
1 & -4-3i\sqrt{6} & -4+3i\sqrt{6} & -4+2i\sqrt{6} & 52+4i\sqrt{6} \\ 
1 & 8+i\sqrt{6} & 8-i\sqrt{6} & -4-3i\sqrt{6} & -14-28i\sqrt{6} \\ 
1 & -4+2i\sqrt{6} & -4-2i\sqrt{6} & -4-3i\sqrt{6} & 52+4i\sqrt{6} \\ 
1 & -4-3i\sqrt{6} & -4+3i\sqrt{6} & -4-3i\sqrt{6} & -38+24i\sqrt{6}%
\end{array}%
\right.  \notag \\
& \left. 
\begin{array}{cccc}
70 & 8-i\sqrt{6} & 70 & 58-16i\sqrt{6} \\ 
-20-20i\sqrt{6} & 8-i\sqrt{6} & -20+20i\sqrt{6} & -44-12i\sqrt{6} \\ 
-50+20i\sqrt{6} & 8-i\sqrt{6} & -50-20i\sqrt{6} & -14+28i\sqrt{6} \\ 
-20+20i\sqrt{6} & -4-2i\sqrt{6} & -20-20i\sqrt{6} & -44-12i\sqrt{6} \\ 
40 & -4-2i\sqrt{6} & 40 & -8+16i\sqrt{6} \\ 
-20-20i\sqrt{6} & -4-2i\sqrt{6} & -20+20i\sqrt{6} & 52-4i\sqrt{6} \\ 
-50-20i\sqrt{6} & -4+3i\sqrt{6} & -50+20i\sqrt{6} & -14+28i\sqrt{6} \\ 
-20+20i\sqrt{6} & -4+3i\sqrt{6} & -20-20i\sqrt{6} & 52-4i\sqrt{6} \\ 
70 & -4+3i\sqrt{6} & 70 & -38-24i\sqrt{6}%
\end{array}%
\right] ,  \label{S9}
\end{align}

\begin{align}
\mathbf{D}_{9}^{A}& =\limfunc{diag}\left[ 36,6i\sqrt{6},-6i\sqrt{6}%
,0,0,0,0,0,0\right] ,  \notag \\
\mathbf{D}_{9}^{B}& =\limfunc{diag}\left[ 36,0,0,6i\sqrt{6},0,0,-6i\sqrt{6}%
,0,0\right] .
\end{align}%
$\allowbreak $Again, the eigenvalues of the regular matrices $\mathbf{A}_{9}$
and $\mathbf{B}_{9}$ are $\mu _{9}$ and $\pm $pairs but in a different order
as related by (\ref{PSP}) with (\ref{P9}). Also, (\ref{PSP}) for $\mathbf{S}%
_{9}$ can be verified. According to (\ref{MJord}) and (\ref{Dn2}), $\mathbf{M%
}_{9}^{A}$ and $\mathbf{M}_{9}^{B}$ have the eigenvector matrix $\mathbf{S}%
_{9}$ and eigenvalues%
\begin{equation}
360,6i\sqrt{6},-6i\sqrt{6},54i\sqrt{6},-54i\sqrt{6},0,0,0,0,  \label{ev}
\end{equation}%
as can be verified directly. These eigenvalues agree with those given by
Rogers, et.\negthinspace\ al \cite{ROGE}.\newpage

Next, we construct the SVD matrices of $\mathbf{A}_{9}$, $\mathbf{B}_{9}$, $%
\mathbf{M}_{9}^{A},$ and $\mathbf{M}_{9}^{B}$ from the following SVD
matrices of $\mathbf{M}_{3}$ and $\mathbf{E}_{3}$:%
\begin{align}
\mathbf{U}_{3}& =\left[ 
\begin{array}{ccc}
\frac{1}{3}\sqrt{3} & \frac{1}{2}\sqrt{2} & \frac{1}{6}\sqrt{6} \\ 
\frac{1}{3}\sqrt{3} & 0 & -\frac{1}{3}\sqrt{6} \\ 
\frac{1}{3}\sqrt{3} & -\frac{1}{2}\sqrt{2} & \frac{1}{6}\sqrt{6}%
\end{array}%
\right] ,  \notag \\
\mathbf{V}_{3}& =\left[ 
\begin{array}{ccc}
\frac{1}{3}\sqrt{3} & -\frac{1}{6}\sqrt{6} & \frac{1}{2}\sqrt{2} \\ 
\frac{1}{3}\sqrt{3} & \frac{1}{3}\sqrt{6} & 0 \\ 
\frac{1}{3}\sqrt{3} & -\frac{1}{6}\sqrt{6} & -\frac{1}{2}\sqrt{2}%
\end{array}%
\right] ,  \label{UV3} \\
\mathbf{\Sigma }_{M3}& =\limfunc{diag}\left[ 12,4\sqrt{3},2\sqrt{3}\right] ,
\notag \\
\mathbf{\Sigma }_{E3}& =\mathbf{D}_{E3}=\limfunc{diag}\left[ 3,0,0\right] . 
\notag
\end{align}%
The SVD matrices for $\mathbf{A}_{9},\ \mathbf{B}_{9},$ $\mathbf{M}_{9}^{A},$
and $\mathbf{M}_{9}^{B}$ are obtained from (\ref{SVDAB}) and (\ref{SigAB}) as%
\begin{equation}
\mathbf{U}_{9}=\frac{1}{6}\left[ 
\begin{array}{ccccccccc}
2 & \sqrt{6} & \sqrt{2} & \sqrt{6} & 3 & \sqrt{3} & \sqrt{2} & \sqrt{3} & 1
\\ 
2 & 0 & -2\sqrt{2} & \sqrt{6} & 0 & -2\sqrt{3} & \sqrt{2} & 0 & -2 \\ 
2 & -\sqrt{6} & \sqrt{2} & \sqrt{6} & -3 & \sqrt{3} & \sqrt{2} & -\sqrt{3} & 
1 \\ 
2 & \sqrt{6} & \sqrt{2} & 0 & 0 & 0 & -2\sqrt{2} & -2\sqrt{3} & -2 \\ 
2 & 0 & -2\sqrt{2} & 0 & 0 & 0 & -2\sqrt{2} & 0 & 4 \\ 
2 & -\sqrt{6} & \sqrt{2} & 0 & 0 & 0 & -2\sqrt{2} & 2\sqrt{3} & -2 \\ 
2 & \sqrt{6} & \sqrt{2} & -\sqrt{6} & -3 & -\sqrt{3} & \sqrt{2} & \sqrt{3} & 
1 \\ 
2 & 0 & -2\sqrt{2} & -\sqrt{6} & 0 & 2\sqrt{3} & \sqrt{2} & 0 & -2 \\ 
2 & -\sqrt{6} & \sqrt{2} & -\sqrt{6} & 3 & -\sqrt{3} & \sqrt{2} & -\sqrt{3}
& 1%
\end{array}%
\right] ,
\end{equation}%
\begin{equation}
\mathbf{V}_{9}=\frac{1}{6}\left[ 
\begin{array}{ccccccccc}
2 & -\sqrt{2} & \sqrt{6} & -\sqrt{2} & 1 & -\sqrt{3} & \sqrt{6} & -\sqrt{3}
& 3 \\ 
2 & 2\sqrt{2} & 0 & -\sqrt{2} & -2 & 0 & \sqrt{6} & 2\sqrt{3} & 0 \\ 
2 & -\sqrt{2} & -\sqrt{6} & -\sqrt{2} & 1 & \sqrt{3} & \sqrt{6} & -\sqrt{3}
& -3 \\ 
2 & -\sqrt{2} & \sqrt{6} & 2\sqrt{2} & -2 & 2\sqrt{3} & 0 & 0 & 0 \\ 
2 & 2\sqrt{2} & 0 & 2\sqrt{2} & 4 & 0 & 0 & 0 & 0 \\ 
2 & -\sqrt{2} & -\sqrt{6} & 2\sqrt{2} & -2 & -2\sqrt{3} & 0 & 0 & 0 \\ 
2 & -\sqrt{2} & \sqrt{6} & -\sqrt{2} & 1 & -\sqrt{3} & -\sqrt{6} & \sqrt{3}
& -3 \\ 
2 & 2\sqrt{2} & 0 & -\sqrt{2} & -2 & 0 & -\sqrt{6} & -2\sqrt{3} & 0 \\ 
2 & -\sqrt{2} & -\sqrt{6} & -\sqrt{2} & 1 & \sqrt{3} & -\sqrt{6} & \sqrt{3}
& 3%
\end{array}%
\right] ,
\end{equation}%
\smallskip 
\begin{align}
\mathbf{\Sigma }_{9}^{A}& =\limfunc{diag}\left[ 36,12\sqrt{3},6\sqrt{3}%
,0,0,0,0,0,0\right] ,  \notag \\
\Sigma _{9}^{B}& =\limfunc{diag}\left[ 36,0,0,12\sqrt{3},0,0,6\sqrt{3},0,0%
\right] ,  \notag \\
\mathbf{\Sigma }_{M9}^{A}& =\limfunc{diag}\left[ 360,12\sqrt{3},6\sqrt{3},108%
\sqrt{3},0,0,54\sqrt{3},0,0\right] , \\
\mathbf{\Sigma }_{M9}^{B}& =\limfunc{diag}\left[ 360,108\sqrt{3},54\sqrt{3}%
,12\sqrt{3},0,0,6\sqrt{3},0,0\right]  \notag
\end{align}%
which can be verified directly from their SVD definitions (\ref{Msvd}).
Also, (\ref{PUP}) can be verified. The singular values for $\mathbf{M}%
_{9}^{A}$ and $\mathbf{M}_{9}^{B}$ agree with those given in \cite{ROGE}%
.\newpage

\noindent \textbf{Order 12.} We start with the regular magic squares used by
Rogers, et.\negthinspace\ al \cite{ROGE}, namely%
\begin{equation}
\mathbf{M}_{n}=\mathbf{M}_{3}=\left[ 
\begin{array}{ccc}
3 & 8 & 1 \\ 
2 & 4 & 6 \\ 
7 & 0 & 5%
\end{array}%
\right] ,\quad \mathbf{M}_{m}=\mathbf{M}_{4}=\left[ 
\begin{array}{cccc}
4 & 3 & 15 & 8 \\ 
10 & 13 & 1 & 6 \\ 
9 & 14 & 2 & 5 \\ 
7 & 0 & 12 & 11%
\end{array}%
\right] .  \label{M3M4}
\end{equation}%
From (\ref{AB}) we form the order-12, commuting regular, unnatural, magic
squares%
\begin{equation}
\mathbf{A}_{12}=\left[ 
\begin{array}{cccccccccccc}
3 & 8 & 1 & 3 & 8 & 1 & 3 & 8 & 1 & 3 & 8 & 1 \\ 
2 & 4 & 6 & 2 & 4 & 6 & 2 & 4 & 6 & 2 & 4 & 6 \\ 
7 & 0 & 5 & 7 & 0 & 5 & 7 & 0 & 5 & 7 & 0 & 5 \\ 
3 & 8 & 1 & 3 & 8 & 1 & 3 & 8 & 1 & 3 & 8 & 1 \\ 
2 & 4 & 6 & 2 & 4 & 6 & 2 & 4 & 6 & 2 & 4 & 6 \\ 
7 & 0 & 5 & 7 & 0 & 5 & 7 & 0 & 5 & 7 & 0 & 5 \\ 
3 & 8 & 1 & 3 & 8 & 1 & 3 & 8 & 1 & 3 & 8 & 1 \\ 
2 & 4 & 6 & 2 & 4 & 6 & 2 & 4 & 6 & 2 & 4 & 6 \\ 
7 & 0 & 5 & 7 & 0 & 5 & 7 & 0 & 5 & 7 & 0 & 5 \\ 
3 & 8 & 1 & 3 & 8 & 1 & 3 & 8 & 1 & 3 & 8 & 1 \\ 
2 & 4 & 6 & 2 & 4 & 6 & 2 & 4 & 6 & 2 & 4 & 6 \\ 
7 & 0 & 5 & 7 & 0 & 5 & 7 & 0 & 5 & 7 & 0 & 5%
\end{array}%
\right] ,
\end{equation}%
\medskip 
\begin{equation}
\mathbf{B}_{12}=\left[ 
\begin{array}{cccccccccccc}
4 & 4 & 4 & 3 & 3 & 3 & 15 & 15 & 15 & 8 & 8 & 8 \\ 
4 & 4 & 4 & 3 & 3 & 3 & 15 & 15 & 15 & 8 & 8 & 8 \\ 
4 & 4 & 4 & 3 & 3 & 3 & 15 & 15 & 15 & 8 & 8 & 8 \\ 
10 & 10 & 10 & 13 & 13 & 13 & 1 & 1 & 1 & 6 & 6 & 6 \\ 
10 & 10 & 10 & 13 & 13 & 13 & 1 & 1 & 1 & 6 & 6 & 6 \\ 
10 & 10 & 10 & 13 & 13 & 13 & 1 & 1 & 1 & 6 & 6 & 6 \\ 
9 & 9 & 9 & 14 & 14 & 14 & 2 & 2 & 2 & 5 & 5 & 5 \\ 
9 & 9 & 9 & 14 & 14 & 14 & 2 & 2 & 2 & 5 & 5 & 5 \\ 
9 & 9 & 9 & 14 & 14 & 14 & 2 & 2 & 2 & 5 & 5 & 5 \\ 
7 & 7 & 7 & 0 & 0 & 0 & 12 & 12 & 12 & 11 & 11 & 11 \\ 
7 & 7 & 7 & 0 & 0 & 0 & 12 & 12 & 12 & 11 & 11 & 11 \\ 
7 & 7 & 7 & 0 & 0 & 0 & 12 & 12 & 12 & 11 & 11 & 11%
\end{array}%
\right] .
\end{equation}%
Then, (\ref{Mn2s}) gives the commuting, regular, magic squares%
\begin{equation}
\mathbf{M}_{12}^{A}=\left[ 
\begin{array}{cccccccccccc}
39 & 44 & 37 & 30 & 35 & 28 & 138 & 143 & 136 & 75 & 80 & 73 \\ 
38 & 40 & 42 & 29 & 31 & 33 & 137 & 139 & 141 & 74 & 76 & 78 \\ 
43 & 36 & 41 & 34 & 27 & 32 & 142 & 135 & 140 & 79 & 72 & 77 \\ 
93 & 98 & 91 & 120 & 125 & 118 & 12 & 17 & 10 & 57 & 62 & 55 \\ 
92 & 94 & 96 & 119 & 121 & 123 & 11 & 13 & 15 & 56 & 58 & 60 \\ 
97 & 90 & 95 & 124 & 117 & 122 & 16 & 9 & 14 & 61 & 54 & 59 \\ 
84 & 89 & 82 & 129 & 134 & 127 & 21 & 26 & 19 & 48 & 53 & 46 \\ 
83 & 85 & 87 & 128 & 130 & 132 & 20 & 22 & 24 & 47 & 49 & 51 \\ 
88 & 81 & 86 & 133 & 126 & 131 & 25 & 18 & 23 & 52 & 45 & 50 \\ 
66 & 71 & 64 & 3 & 8 & 1 & 111 & 116 & 109 & 102 & 107 & 100 \\ 
65 & 67 & 69 & 2 & 4 & 6 & 110 & 112 & 114 & 101 & 103 & 105 \\ 
70 & 63 & 68 & 7 & 0 & 5 & 115 & 108 & 113 & 106 & 99 & 104%
\end{array}%
\right] ,
\end{equation}%
\begin{equation}
\mathbf{M}_{12}^{B}=\left[ 
\begin{array}{cccccccccccc}
52 & 132 & 20 & 51 & 131 & 19 & 63 & 143 & 31 & 56 & 136 & 24 \\ 
36 & 68 & 100 & 35 & 67 & 99 & 47 & 79 & 111 & 40 & 72 & 104 \\ 
116 & 4 & 84 & 115 & 3 & 83 & 127 & 15 & 95 & 120 & 8 & 88 \\ 
58 & 138 & 26 & 61 & 141 & 29 & 49 & 129 & 17 & 54 & 134 & 22 \\ 
42 & 74 & 106 & 45 & 77 & 109 & 33 & 65 & 97 & 38 & 70 & 102 \\ 
122 & 10 & 90 & 125 & 13 & 93 & 113 & 1 & 81 & 118 & 6 & 86 \\ 
57 & 137 & 25 & 62 & 142 & 30 & 50 & 130 & 18 & 53 & 133 & 21 \\ 
41 & 73 & 105 & 46 & 78 & 110 & 34 & 66 & 98 & 37 & 69 & 101 \\ 
121 & 9 & 89 & 126 & 14 & 94 & 114 & 2 & 82 & 117 & 5 & 85 \\ 
55 & 135 & 23 & 48 & 128 & 16 & 60 & 140 & 28 & 59 & 139 & 27 \\ 
39 & 71 & 103 & 32 & 64 & 96 & 44 & 76 & 108 & 43 & 75 & 107 \\ 
119 & 7 & 87 & 112 & 0 & 80 & 124 & 12 & 92 & 123 & 11 & 91%
\end{array}%
\right] .
\end{equation}

The Jordan-form matrices of $\mathbf{M}_{3}$ and $\mathbf{E}_{3}$ are given
by (\ref{Jord3}) and those of $\mathbf{M}_{4}$ are%
\begin{equation}
\mathbf{S}_{4}=\left[ 
\begin{array}{cccc}
1 & 48 & -14 & 3 \\ 
1 & -16 & 10 & -1 \\ 
1 & 16 & 6 & -1 \\ 
1 & -48 & -2 & -1%
\end{array}%
\right] ,\quad \mathbf{J}_{4}=\left[ 
\begin{array}{cccc}
30 & 0 & 0 & 0 \\ 
0 & 0 & 1 & 0 \\ 
0 & 0 & 0 & 1 \\ 
0 & 0 & 0 & 0%
\end{array}%
\right] ,
\end{equation}%
where $\mathbf{S}_{4}$ has a generalized eigenvector.$\mathbf{\ }$Then, (\ref%
{SmnD}) results in%
\begin{align}
\mathbf{S}_{12}& =\left[ 
\begin{array}{cccccc}
1 & 8+i\sqrt{6} & 8-i\sqrt{6} & 48 & 384+48i\sqrt{6} & 384-48i\sqrt{6} \\ 
1 & -4+2i\sqrt{6} & -4-2i\sqrt{6} & 48 & -192+96i\sqrt{6} & -192-96i\sqrt{6}
\\ 
1 & -4-3i\sqrt{6} & -4+3i\sqrt{6} & 48 & -192-144i\sqrt{6} & -192+144i\sqrt{6%
} \\ 
1 & 8+i\sqrt{6} & 8-i\sqrt{6} & -16 & -128-16i\sqrt{6} & -128+16i\sqrt{6} \\ 
1 & -4+2i\sqrt{6} & -4-2i\sqrt{6} & -16 & 64-32i\sqrt{6} & 64+32i\sqrt{6} \\ 
1 & -4-3i\sqrt{6} & -4+3i\sqrt{6} & -16 & 64+48i\sqrt{6} & 64-48i\sqrt{6} \\ 
1 & 8+i\sqrt{6} & 8-i\sqrt{6} & 16 & 128+16i\sqrt{6} & 128-16i\sqrt{6} \\ 
1 & -4+2i\sqrt{6} & -4-2i\sqrt{6} & 16 & -64+32i\sqrt{6} & -64-32i\sqrt{6}
\\ 
1 & -4-3i\sqrt{6} & -4+3i\sqrt{6} & 16 & -64-48i\sqrt{6} & -64+48i\sqrt{6}
\\ 
1 & 8+i\sqrt{6} & 8-i\sqrt{6} & -48 & -384-48i\sqrt{6} & -384+48i\sqrt{6} \\ 
1 & -4+2i\sqrt{6} & -4-2i\sqrt{6} & -48 & 192-96i\sqrt{6} & 192+96i\sqrt{6}
\\ 
1 & -4-3i\sqrt{6} & -4+3i\sqrt{6} & -48 & 192+144i\sqrt{6} & 192-144i\sqrt{6}%
\end{array}%
\right.  \notag \\
& \left. 
\begin{array}{cccccc}
-14 & -112-14i\sqrt{6} & -112+14i\sqrt{6} & 3 & 24+3i\sqrt{6} & 24-3i\sqrt{6}
\\ 
-14 & 56-28i\sqrt{6} & 56+28i\sqrt{6} & 3 & -12+6i\sqrt{6} & -12-6i\sqrt{6}
\\ 
-14 & 56+42i\sqrt{6} & 56-42i\sqrt{6} & 3 & -12-9i\sqrt{6} & -12+9i\sqrt{6}
\\ 
10 & 80+10i\sqrt{6} & 80-10i\sqrt{6} & -1 & -8-i\sqrt{6} & -8+i\sqrt{6} \\ 
10 & -40+20i\sqrt{6} & -40-20i\sqrt{6} & -1 & 4-2i\sqrt{6} & 4+2i\sqrt{6} \\ 
10 & -40-30i\sqrt{6} & -40+30i\sqrt{6} & -1 & 4+3i\sqrt{6} & 4-3i\sqrt{6} \\ 
6 & 48+6i\sqrt{6} & 48-6i\sqrt{6} & -1 & -8-i\sqrt{6} & -8+i\sqrt{6} \\ 
6 & -24+12i\sqrt{6} & -24-12i\sqrt{6} & -1 & 4-2i\sqrt{6} & 4+2i\sqrt{6} \\ 
6 & -24-18i\sqrt{6} & -24+18i\sqrt{6} & -1 & 4+3i\sqrt{6} & 4-3i\sqrt{6} \\ 
-2 & -16-2i\sqrt{6} & -16+2i\sqrt{6} & -1 & -8-i\sqrt{6} & -8+i\sqrt{6} \\ 
-2 & 8-4i\sqrt{6} & 8+4i\sqrt{6} & -1 & 4-2i\sqrt{6} & 4+2i\sqrt{6} \\ 
-2 & 8+6i\sqrt{6} & 8-6i\sqrt{6} & -1 & 4+3i\sqrt{6} & 4-3i\sqrt{6}%
\end{array}%
\right] ,
\end{align}%
\begin{align}
\mathbf{J}_{A12}& =\limfunc{diag}\left[ 12,8i\sqrt{6},-8i\sqrt{6}.0,\ldots ,0%
\right] ,  \notag \\
\mathbf{J}_{B12}& =\left[ 
\begin{array}{cccccccccccc}
90 & 0 & 0 & 0 & 0 & 0 & 0 & 0 & 0 & 0 & 0 & 0 \\ 
0 & 0 & 0 & 0 & 0 & 0 & 0 & 0 & 0 & 0 & 0 & 0 \\ 
0 & 0 & 0 & 0 & 0 & 0 & 0 & 0 & 0 & 0 & 0 & 0 \\ 
0 & 0 & 0 & 0 & 0 & 0 & 3 & 0 & 0 & 0 & 0 & 0 \\ 
0 & 0 & 0 & 0 & 0 & 0 & 0 & 0 & 0 & 0 & 0 & 0 \\ 
0 & 0 & 0 & 0 & 0 & 0 & 0 & 0 & 0 & 0 & 0 & 0 \\ 
0 & 0 & 0 & 0 & 0 & 0 & 0 & 0 & 0 & 3 & 0 & 0 \\ 
0 & 0 & 0 & 0 & 0 & 0 & 0 & 0 & 0 & 0 & 0 & 0 \\ 
0 & 0 & 0 & 0 & 0 & 0 & 0 & 0 & 0 & 0 & 0 & 0 \\ 
0 & 0 & 0 & 0 & 0 & 0 & 0 & 0 & 0 & 0 & 0 & 0 \\ 
0 & 0 & 0 & 0 & 0 & 0 & 0 & 0 & 0 & 0 & 0 & 0 \\ 
0 & 0 & 0 & 0 & 0 & 0 & 0 & 0 & 0 & 0 & 0 & 0%
\end{array}%
\right] .
\end{align}%
Although $\mathbf{J}_{B12}$ is not in standard form, it can be brought there
by rearranging and scaling the generalized eigenvectors in $\mathbf{S}_{12}$
without affecting the eigenvalues in $\mathbf{J}_{A12}$. From (\ref{Dn2}),
the nonzero eigenvalues of $\mathbf{M}_{12}^{A}$ and $\mathbf{M}_{12}^{B}$
are found to be%
\begin{equation}
\mathbf{M}_{12}^{A}:858,8i\sqrt{6},-8i\sqrt{6},\quad \mathbf{M}%
_{12}^{B}:858,128i\sqrt{6},-128i\sqrt{6}.
\end{equation}%
The nonzero eigenvalues of $\mathbf{M}_{12}^{A}$ agree with those given by
Rogers, et.\negthinspace\ al \cite{ROGE} who do not construct $\mathbf{M}%
_{12}^{B}.$

The SVD matrices of $\mathbf{M}_{3}$ are given by (\ref{UV3}) and those of $%
\mathbf{M}_{4}$ are%
\begin{align}
\mathbf{U}_{4}& =\frac{1}{10}\left[ 
\begin{array}{cccc}
5 & -5 & 3\sqrt{5} & \sqrt{5} \\ 
5 & 5 & -\sqrt{5} & 3\sqrt{5} \\ 
5 & 5 & \sqrt{5} & -3\sqrt{5} \\ 
5 & -5 & -3\sqrt{5} & -\sqrt{5}%
\end{array}%
\right] ,  \notag \\
\mathbf{V}_{4}& =\frac{1}{10}\left[ 
\begin{array}{cccc}
5 & \sqrt{5} & -5 & -3\sqrt{5} \\ 
5 & 3\sqrt{5} & 5 & \sqrt{5} \\ 
5 & -3\sqrt{5} & 5 & -\sqrt{5} \\ 
5 & -\sqrt{5} & -5 & 3\sqrt{5}%
\end{array}%
\right] ,  \label{UV4} \\
\mathbf{\Sigma }_{4}& =\limfunc{diag}\left[ 30,8\sqrt{5},2\sqrt{5},0\right] .
\notag
\end{align}%
$\allowbreak $The $\mathbf{U}_{12}$ and $\mathbf{V}_{12}$ matrices for $%
\mathbf{A}_{12}$, $\mathbf{B}_{12}$, $\mathbf{M}_{12}^{A}$, and $\mathbf{M}%
_{12}^{B}$ can be obtained from (\ref{SVDAB}) with (\ref{UV3}) and (\ref{UV4}%
). The singular value matrices from (\ref{SVDAB}), (\ref{UV3}), (\ref{UV4})
and (\ref{SigAB}) are\footnote{%
The singular values obtained for $\mathbf{M}_{12}^{\left( A\right) }$
(verified by MAPLE).agree with those of Rogers, et. \negthinspace
\negthinspace al \cite{ROGE} except for their $24\sqrt{5}$ and $6\sqrt{5}$
instead of our $216\sqrt{5}$ and $54\sqrt{5}$.}%
\begin{align}
\mathbf{\Sigma }_{12}^{A}& =\limfunc{diag}\left[ 48,16\sqrt{3},8\sqrt{3}%
,0,0,...,0\right] ,  \notag \\
\mathbf{\Sigma }_{12}^{B}& =\limfunc{diag}\left[ 90,0,0,24\sqrt{5},0,0,6%
\sqrt{5},0,0,...,0\right] ,  \notag \\
\mathbf{\Sigma }_{M12}^{A}& =\limfunc{diag}\left[ 858,16\sqrt{3},8\sqrt{3}%
,216\sqrt{5},0,0,54\sqrt{5},0,0,\ldots ,0\right] , \\
\mathbf{\Sigma }_{M12}^{B}& =\limfunc{diag}\left[ 858,256\sqrt{3},128\sqrt{3}%
,24\sqrt{5},0,0,6\sqrt{5},0,0,\ldots ,0\right] .  \notag
\end{align}%
\newpage 

In addition, as noted by Rogers, et. \negthinspace al \cite{ROGE}, a
different pair of $\mathbf{M}_{12}^{A}$ and $\mathbf{M}_{12}^{B}$ can be
constructed by interchanging $\mathbf{M}_{m}$ and $\mathbf{M}_{n}$ in (\ref%
{M3M4}), i.e.%
\begin{equation}
\mathbf{M}_{n}=\mathbf{M}_{4}=\left[ 
\begin{array}{cccc}
4 & 3 & 15 & 8 \\ 
10 & 13 & 1 & 6 \\ 
9 & 14 & 2 & 5 \\ 
7 & 0 & 12 & 11%
\end{array}%
\right] ,\quad \mathbf{M}_{m}=\mathbf{M}_{3}=\left[ 
\begin{array}{ccc}
3 & 8 & 1 \\ 
2 & 4 & 6 \\ 
7 & 0 & 5%
\end{array}%
\right] ,
\end{equation}%
leading to the commuting pair of regular magic squares%
\begin{equation}
\mathbf{\hat{M}}_{12}^{A}=\left[ 
\begin{array}{cccccccccccc}
52 & 51 & 63 & 56 & 132 & 131 & 143 & 136 & 20 & 19 & 31 & 24 \\ 
58 & 61 & 49 & 54 & 138 & 141 & 129 & 134 & 26 & 29 & 17 & 22 \\ 
57 & 62 & 50 & 53 & 137 & 142 & 130 & 133 & 25 & 30 & 18 & 21 \\ 
55 & 48 & 60 & 59 & 135 & 128 & 140 & 139 & 23 & 16 & 28 & 27 \\ 
36 & 35 & 47 & 40 & 68 & 67 & 79 & 72 & 100 & 99 & 111 & 104 \\ 
42 & 45 & 33 & 38 & 74 & 77 & 65 & 70 & 106 & 109 & 97 & 102 \\ 
41 & 46 & 34 & 37 & 73 & 78 & 66 & 69 & 105 & 110 & 98 & 101 \\ 
39 & 32 & 44 & 43 & 71 & 64 & 76 & 75 & 103 & 96 & 108 & 107 \\ 
116 & 115 & 127 & 120 & 4 & 3 & 15 & 8 & 84 & 83 & 95 & 88 \\ 
122 & 125 & 113 & 118 & 10 & 13 & 1 & 6 & 90 & 93 & 81 & 86 \\ 
121 & 126 & 114 & 117 & 9 & 14 & 2 & 5 & 89 & 94 & 82 & 85 \\ 
119 & 112 & 124 & 123 & 7 & 0 & 12 & 11 & 87 & 80 & 92 & 91%
\end{array}%
\right] ,
\end{equation}%
\medskip 
\begin{equation}
\mathbf{\hat{M}}_{12}^{B}=\left[ 
\begin{array}{cccccccccccc}
39 & 30 & 138 & 75 & 44 & 35 & 143 & 80 & 37 & 28 & 136 & 73 \\ 
93 & 120 & 12 & 57 & 98 & 125 & 17 & 62 & 91 & 118 & 10 & 55 \\ 
84 & 129 & 21 & 48 & 89 & 134 & 26 & 53 & 82 & 127 & 19 & 46 \\ 
66 & 3 & 111 & 102 & 71 & 8 & 116 & 107 & 64 & 1 & 109 & 100 \\ 
38 & 29 & 137 & 74 & 40 & 31 & 139 & 76 & 42 & 33 & 141 & 78 \\ 
92 & 119 & 11 & 56 & 94 & 121 & 13 & 58 & 96 & 123 & 15 & 60 \\ 
83 & 128 & 20 & 47 & 85 & 130 & 22 & 49 & 87 & 132 & 24 & 51 \\ 
65 & 2 & 110 & 101 & 67 & 4 & 112 & 103 & 69 & 6 & 114 & 105 \\ 
43 & 34 & 142 & 79 & 36 & 27 & 135 & 72 & 41 & 32 & 140 & 77 \\ 
97 & 124 & 16 & 61 & 90 & 117 & 9 & 54 & 95 & 122 & 14 & 59 \\ 
88 & 133 & 25 & 52 & 81 & 126 & 18 & 45 & 86 & 131 & 23 & 50 \\ 
70 & 7 & 115 & 106 & 63 & 0 & 108 & 99 & 68 & 5 & 113 & 104%
\end{array}%
\right] .
\end{equation}%
The matrices in the Jordan form and SVD of $\mathbf{\hat{M}}_{12}^{A}$ and $%
\mathbf{\hat{M}}_{12}^{B}$ can be constructed as before. It follows from
formulas in Sections 3 and 4 that the eigenvalues and singular values for $%
\mathbf{\hat{M}}_{12}^{A}$ are the same as those for $\mathbf{M}_{12}^{B}$
and the eigenvalues and singular values for $\mathbf{\hat{M}}_{12}^{B}$ are
the same as those for $\mathbf{M}_{12}^{A}$ (verified by MAPLE). However,
their respective eigenvalue matrices are different and they do not commute.

Again, as noted in Section 2, a huge number of pairs of commuting magic
squares $\mathbf{M}_{12}^{A}$ and $\mathbf{M}_{12}^{B}$ can be constructed
according to (\ref{Mn2s}) using various combinations of phases of $\mathbf{M}%
_{3}$ as subsquares in generalized $\mathbf{\tilde{A}}_{12}$ and any $%
\mathbf{M}_{4}$ as the basis for $\mathbf{B}_{12}$ in (\ref{AB}). This same
construction applies to $\mathbf{\hat{M}}_{12}^{A}$ and $\mathbf{\hat{M}}%
_{12}^{B}$ using various combinations of phases of $\mathbf{M}_{4}$ as
subsquares in generalized $\mathbf{\tilde{A}}_{12}$ and any $\mathbf{M}_{3}$
as the basis for $\mathbf{B}_{12}$ in (\ref{AB}). \newpage

\noindent \textbf{Order 16.} We start with the order-$4,$ pandiagonal, magic
square equivalent to the one used by Chan and Loly \cite{CHAN}, namely%
\begin{equation}
\mathbf{M}_{4}=\left[ 
\begin{array}{cccc}
13 & 6 & 11 & 0 \\ 
10 & 1 & 12 & 7 \\ 
4 & 15 & 2 & 9 \\ 
3 & 8 & 5 & 14%
\end{array}%
\right] .
\end{equation}%
By compounding according to (\ref{AB}) with $\mathbf{M}_{m}=\mathbf{M}_{n}=%
\mathbf{M}_{4},$ we form the following order-16, commuting, pandiagonal,
unnatural magic squares%
\begin{equation}
\mathbf{A}_{16}=%
\begin{tabular}{|cccccccccccccccc|}
\hline
{\small 13} & \textbf{6} & {\small 11} & {\small 0} & {\small 13} & {\small 6%
} & {\small 11} & {\small 0} & {\small 13} & \textbf{6} & {\small 11} & 
{\small 0} & {\small 13} & {\small 6} & {\small 11} & {\small 0} \\ 
{\small 10} & {\small 1} & \textbf{12} & {\small 7} & {\small 10} & {\small 1%
} & {\small 12} & {\small 7} & \textbf{10} & {\small 1} & {\small 12} & 
{\small 7} & {\small 10} & {\small 1} & {\small 12} & {\small 7} \\ 
{\small 4} & {\small 15} & {\small 2} & \textbf{9} & {\small 4} & {\small 15}
& {\small 2} & \textbf{9} & {\small 4} & {\small 15} & {\small 2} & {\small 9%
} & {\small 4} & {\small 15} & {\small 2} & {\small 9} \\ 
{\small 3} & {\small 8} & {\small 5} & {\small 14} & \textbf{3} & {\small 8}
& \textbf{5} & {\small 14} & {\small 3} & {\small 8} & {\small 5} & {\small %
14} & {\small 3} & {\small 8} & {\small 5} & {\small 14} \\ 
{\small 13} & {\small 6} & {\small 11} & {\small 0} & {\small 13} & \textbf{6%
} & {\small 11} & {\small 0} & {\small 13} & {\small 6} & {\small 11} & 
{\small 0} & {\small 13} & {\small 6} & {\small 11} & {\small 0} \\ 
{\small 10} & {\small 1} & {\small 12} & {\small 7} & \textbf{10} & {\small 1%
} & \textbf{12} & {\small 7} & {\small 10} & {\small 1} & {\small 12} & 
{\small 7} & {\small 10} & {\small 1} & {\small 12} & {\small 7} \\ 
{\small 4} & {\small 15} & {\small 2} & \textbf{9} & {\small 4} & {\small 15}
& {\small 2} & \textbf{9} & {\small 4} & {\small 15} & {\small 2} & {\small 9%
} & {\small 4} & {\small 15} & {\small 2} & {\small 9} \\ 
{\small 3} & {\small 8} & \textbf{5} & {\small 14} & {\small 3} & {\small 8}
& {\small 5} & {\small 14} & \textbf{3} & {\small 8} & {\small 5} & {\small %
14} & {\small 3} & {\small 8} & {\small 5} & {\small 14} \\ 
{\small 13} & \textbf{6} & {\small 11} & {\small 0} & {\small 13} & {\small 6%
} & {\small 11} & {\small 0} & {\small 13} & \textbf{6} & {\small 11} & 
{\small 0} & {\small 13} & {\small 6} & {\small 11} & {\small 0} \\ 
\textbf{10} & {\small 1} & {\small 12} & {\small 7} & {\small 10} & {\small 1%
} & {\small 12} & {\small 7} & {\small 10} & {\small 1} & \textbf{12} & 
{\small 7} & {\small 10} & {\small 1} & {\small 12} & {\small 7} \\ 
{\small 4} & {\small 15} & {\small 2} & {\small 9} & {\small 4} & {\small 15}
& {\small 2} & {\small 9} & {\small 4} & {\small 15} & {\small 2} & \textbf{9%
} & {\small 4} & {\small 15} & {\small 2} & \textbf{9} \\ 
{\small 3} & {\small 8} & {\small 5} & {\small 14} & {\small 3} & {\small 8}
& {\small 5} & {\small 14} & {\small 3} & {\small 8} & {\small 5} & {\small %
14} & \textbf{3} & {\small 8} & \textbf{5} & {\small 14} \\ 
{\small 13} & {\small 6} & {\small 11} & {\small 0} & {\small 13} & {\small 6%
} & {\small 11} & {\small 0} & {\small 13} & {\small 6} & {\small 11} & 
{\small 0} & {\small 13} & \textbf{6} & {\small 11} & {\small 0} \\ 
{\small 10} & {\small 1} & {\small 12} & {\small 7} & {\small 10} & {\small 1%
} & {\small 12} & {\small 7} & {\small 10} & {\small 1} & {\small 12} & 
{\small 7} & \textbf{10} & {\small 1} & \textbf{12} & {\small 7} \\ 
{\small 4} & {\small 15} & {\small 2} & {\small 9} & {\small 4} & {\small 15}
& {\small 2} & {\small 9} & {\small 4} & {\small 15} & {\small 2} & \textbf{9%
} & {\small 4} & {\small 15} & {\small 2} & \textbf{9} \\ 
\textbf{3} & {\small 8} & {\small 5} & {\small 14} & {\small 3} & {\small 8}
& {\small 5} & {\small 14} & {\small 3} & {\small 8} & \textbf{5} & {\small %
14} & {\small 3} & {\small 8} & {\small 5} & {\small 14} \\ \hline
\end{tabular}%
\text{\thinspace ,}
\end{equation}%
\medskip 
\begin{equation}
\mathbf{B}_{16}=%
\begin{tabular}{|cccccccccccccccc|}
\hline
{\small 13} & {\small 13} & \textbf{13} & {\small 13} & {\small 6} & {\small %
6} & {\small 6} & {\small 6} & {\small 11} & {\small 11} & \textbf{11} & 
{\small 11} & {\small 0} & {\small 0} & {\small 0} & {\small 0} \\ 
{\small 13} & {\small 13} & {\small 13} & \textbf{13} & {\small 6} & {\small %
6} & {\small 6} & {\small 6} & {\small 11} & \textbf{11} & {\small 11} & 
{\small 11} & {\small 0} & {\small 0} & {\small 0} & {\small 0} \\ 
{\small 13} & {\small 13} & {\small 13} & {\small 13} & \textbf{6} & {\small %
6} & {\small 6} & {\small 6} & \textbf{11} & {\small 11} & {\small 11} & 
{\small 11} & {\small 0} & {\small 0} & {\small 0} & {\small 0} \\ 
{\small 13} & {\small 13} & {\small 13} & {\small 13} & {\small 6} & \textbf{%
6} & {\small 6} & \textbf{6} & {\small 11} & {\small 11} & {\small 11} & 
{\small 11} & {\small 0} & {\small 0} & {\small 0} & {\small 0} \\ 
{\small 10} & {\small 10} & {\small 10} & {\small 10} & {\small 1} & {\small %
1} & \textbf{1} & {\small 1} & {\small 12} & {\small 12} & {\small 12} & 
{\small 12} & {\small 7} & {\small 7} & {\small 7} & {\small 7} \\ 
{\small 10} & {\small 10} & {\small 10} & {\small 10} & {\small 1} & \textbf{%
1} & {\small 1} & \textbf{1} & {\small 12} & {\small 12} & {\small 12} & 
{\small 12} & {\small 7} & {\small 7} & {\small 7} & {\small 7} \\ 
{\small 10} & {\small 10} & {\small 10} & {\small 10} & \textbf{1} & {\small %
1} & {\small 1} & {\small 1} & \textbf{12} & {\small 12} & {\small 12} & 
{\small 12} & {\small 7} & {\small 7} & {\small 7} & {\small 7} \\ 
{\small 10} & {\small 10} & {\small 10} & \textbf{10} & {\small 1} & {\small %
1} & {\small 1} & {\small 1} & {\small 12} & \textbf{12} & {\small 12} & 
{\small 12} & {\small 7} & {\small 7} & {\small 7} & {\small 7} \\ 
{\small 4} & {\small 4} & \textbf{4} & {\small 4} & {\small 15} & {\small 15}
& {\small 15} & {\small 15} & {\small 2} & {\small 2} & \textbf{2} & {\small %
2} & {\small 9} & {\small 9} & {\small 9} & {\small 9} \\ 
{\small 4} & \textbf{4} & {\small 4} & {\small 4} & {\small 15} & {\small 15}
& {\small 15} & {\small 15} & {\small 2} & {\small 2} & {\small 2} & \textbf{%
2} & {\small 9} & {\small 9} & {\small 9} & {\small 9} \\ 
\textbf{4} & {\small 4} & {\small 4} & {\small 4} & {\small 15} & {\small 15}
& {\small 15} & {\small 15} & {\small 2} & {\small 2} & {\small 2} & {\small %
2} & \textbf{9} & {\small 9} & {\small 9} & {\small 9} \\ 
{\small 4} & {\small 4} & {\small 4} & {\small 4} & {\small 15} & {\small 15}
& {\small 15} & {\small 15} & {\small 2} & {\small 2} & {\small 2} & {\small %
2} & {\small 9} & \textbf{9} & {\small 9} & \textbf{9} \\ 
{\small 3} & {\small 3} & {\small 3} & {\small 3} & {\small 8} & {\small 8}
& {\small 8} & {\small 8} & {\small 5} & {\small 5} & {\small 5} & {\small 5}
& {\small 14} & {\small 14} & \textbf{14} & {\small 14} \\ 
{\small 3} & {\small 3} & {\small 3} & {\small 3} & {\small 8} & {\small 8}
& {\small 8} & {\small 8} & {\small 5} & {\small 5} & {\small 5} & {\small 5}
& {\small 14} & \textbf{14} & {\small 14} & \textbf{14} \\ 
\textbf{3} & {\small 3} & {\small 3} & {\small 3} & {\small 8} & {\small 8}
& {\small 8} & {\small 8} & {\small 5} & {\small 5} & {\small 5} & {\small 5}
& \textbf{14} & {\small 14} & {\small 14} & {\small 14} \\ 
{\small 3} & \textbf{3} & {\small 3} & {\small 3} & {\small 8} & {\small 8}
& {\small 8} & {\small 8} & {\small 5} & {\small 5} & {\small 5} & \textbf{5}
& {\small 14} & {\small 14} & {\small 14} & {\small 14} \\ \hline
\end{tabular}%
\text{\thinspace .}
\end{equation}%
It is not difficult to verify that these two matrices are pandiagonal by
comparing their diagonals with those of $\mathbf{M}_{4}$ as done by
Eggermont \cite{EGGE}. Two such diagonals of $\mathbf{A}_{16}$ and $\mathbf{B%
}_{16}$ are shown in bold. The elements on the diagonals of $\mathbf{A}_{16}$
are simply four replications of the diagonals of $\mathbf{M}_{4}$. The
elements on the diagonals of $\mathbf{B}_{16}$ are weighted combinations of
two adjacent diagonals of $\mathbf{M}_{4}$. It should be clear that a
similar argument applies to other cases of compounded pandiagonal magic
squares\footnote{%
A lenthy formal proof of the pandiagonality of general $\mathbf{A}_{mn}$ and 
$\mathbf{B}_{mn}$ is possible but it is not given here.}.

From $\mathbf{A}_{16}$ and $\mathbf{B}_{16}$, (\ref{Mn2s}) gives the
commuting, pandiagonal, natural, magic squares%
\begin{gather}
\mathbf{M}_{16}^{A}=  \notag \\
\left[ 
\begin{array}{cccccccccccccccc}
221 & 214 & 219 & 208 & 109 & 102 & 107 & 96 & 189 & 182 & 187 & 176 & 13 & 6
& 11 & 0 \\ 
218 & 209 & 220 & 215 & 106 & 97 & 108 & 103 & 186 & 177 & 188 & 183 & 10 & 1
& 12 & 7 \\ 
212 & 223 & 210 & 217 & 100 & 111 & 98 & 105 & 180 & 191 & 178 & 185 & 4 & 15
& 2 & 9 \\ 
211 & 216 & 213 & 222 & 99 & 104 & 101 & 110 & 179 & 184 & 181 & 190 & 3 & 8
& 5 & 14 \\ 
173 & 166 & 171 & 160 & 29 & 22 & 27 & 16 & 205 & 198 & 203 & 192 & 125 & 118
& 123 & 112 \\ 
170 & 161 & 172 & 167 & 26 & 17 & 28 & 23 & 202 & 193 & 204 & 199 & 122 & 113
& 124 & 119 \\ 
164 & 175 & 162 & 169 & 20 & 31 & 18 & 25 & 196 & 207 & 194 & 201 & 116 & 127
& 114 & 121 \\ 
163 & 168 & 165 & 174 & 19 & 24 & 21 & 30 & 195 & 200 & 197 & 206 & 115 & 120
& 117 & 126 \\ 
77 & 70 & 75 & 64 & 253 & 246 & 251 & 240 & 45 & 38 & 43 & 32 & 157 & 150 & 
155 & 144 \\ 
74 & 65 & 76 & 71 & 250 & 241 & 252 & 247 & 42 & 33 & 44 & 39 & 154 & 145 & 
156 & 151 \\ 
68 & 79 & 66 & 73 & 244 & 255 & 242 & 249 & 36 & 47 & 34 & 41 & 148 & 159 & 
146 & 153 \\ 
67 & 72 & 69 & 78 & 243 & 248 & 245 & 254 & 35 & 40 & 37 & 46 & 147 & 152 & 
149 & 158 \\ 
61 & 54 & 59 & 48 & 141 & 134 & 139 & 128 & 93 & 86 & 91 & 80 & 237 & 230 & 
235 & 224 \\ 
58 & 49 & 60 & 55 & 138 & 129 & 140 & 135 & 90 & 81 & 92 & 87 & 234 & 225 & 
236 & 231 \\ 
52 & 63 & 50 & 57 & 132 & 143 & 130 & 137 & 84 & 95 & 82 & 89 & 228 & 239 & 
226 & 233 \\ 
51 & 56 & 53 & 62 & 131 & 136 & 133 & 142 & 83 & 88 & 85 & 94 & 227 & 232 & 
229 & 238%
\end{array}%
\right] ,\medskip \\
\mathbf{M}_{16}^{B}=  \notag \\
\left[ 
\begin{array}{cccccccccccccccc}
221 & 109 & 189 & 13 & 214 & 102 & 182 & 6 & 219 & 107 & 187 & 11 & 208 & 96
& 176 & 0 \\ 
173 & 29 & 205 & 125 & 166 & 22 & 198 & 118 & 171 & 27 & 203 & 123 & 160 & 16
& 192 & 112 \\ 
77 & 253 & 45 & 157 & 70 & 246 & 38 & 150 & 75 & 251 & 43 & 155 & 64 & 240 & 
32 & 144 \\ 
61 & 141 & 93 & 237 & 54 & 134 & 86 & 230 & 59 & 139 & 91 & 235 & 48 & 128 & 
80 & 224 \\ 
218 & 106 & 186 & 10 & 209 & 97 & 177 & 1 & 220 & 108 & 188 & 12 & 215 & 103
& 183 & 7 \\ 
170 & 26 & 202 & 122 & 161 & 17 & 193 & 113 & 172 & 28 & 204 & 124 & 167 & 23
& 199 & 119 \\ 
74 & 250 & 42 & 154 & 65 & 241 & 33 & 145 & 76 & 252 & 44 & 156 & 71 & 247 & 
39 & 151 \\ 
58 & 138 & 90 & 234 & 49 & 129 & 81 & 225 & 60 & 140 & 92 & 236 & 55 & 135 & 
87 & 231 \\ 
212 & 100 & 180 & 4 & 223 & 111 & 191 & 15 & 210 & 98 & 178 & 2 & 217 & 105
& 185 & 9 \\ 
164 & 20 & 196 & 116 & 175 & 31 & 207 & 127 & 162 & 18 & 194 & 114 & 169 & 25
& 201 & 121 \\ 
68 & 244 & 36 & 148 & 79 & 255 & 47 & 159 & 66 & 242 & 34 & 146 & 73 & 249 & 
41 & 153 \\ 
52 & 132 & 84 & 228 & 63 & 143 & 95 & 239 & 50 & 130 & 82 & 226 & 57 & 137 & 
89 & 233 \\ 
211 & 99 & 179 & 3 & 216 & 104 & 184 & 8 & 213 & 101 & 181 & 5 & 222 & 110 & 
190 & 14 \\ 
163 & 19 & 195 & 115 & 168 & 24 & 200 & 120 & 165 & 21 & 197 & 117 & 174 & 30
& 206 & 126 \\ 
67 & 243 & 35 & 147 & 72 & 248 & 40 & 152 & 69 & 245 & 37 & 149 & 78 & 254 & 
46 & 158 \\ 
51 & 131 & 83 & 227 & 56 & 136 & 88 & 232 & 53 & 133 & 85 & 229 & 62 & 142 & 
94 & 238%
\end{array}%
\right] ,
\end{gather}%
where pandiagonality follows from (\ref{Mn2s}) and the pandiagonality of $%
\mathbf{A}_{16}$ and $\mathbf{B}_{16}.$ The matrix $\mathbf{M}_{16}^{A}$
essentially agrees with the matrix constructed in \cite{CHAN}. The similar
matrices $\mathbf{M}_{16}^{A}$ and $\mathbf{M}_{16}^{B}$ are related by (\ref%
{PermM}) with $\mathbf{P}_{16}$ constructed from (\ref{Pn2}). We leave to
the reader the pleasure (or pain) of constructing the Jordan form and SVD of
these order-16 matrices from our formulas for them.

Again, as noted in Section 2, a huge number of pairs of commuting magic
squares $\mathbf{\tilde{M}}_{16}^{A}$ and $\mathbf{\tilde{M}}_{16}^{B}$ can
be constructed according to (\ref{Mn2s}) using various 16 combinations of
the known 880 $\mathbf{M}_{4}$ magic squares \cite{PICK} as the 16
subsquares of generalized $\mathbf{\tilde{A}}_{16}$ and any $\mathbf{M}_{4}$
as the basis for $\mathbf{B}_{16}$ in (\ref{AB}). Pandiagonal $\mathbf{%
\tilde{M}}_{16}^{A}$ and $\mathbf{\tilde{M}}_{16}^{B}$ of this form are not
possible in general except by first constructing regular $\mathbf{\tilde{M}}%
_{16}^{A}$ and $\mathbf{\tilde{M}}_{16}^{B}$ using any eight regular $%
\mathbf{M}_{4}$ as subsquares of $\mathbf{\tilde{A}}_{16}$ arranged in a
regular block pattern and transforming these $\mathbf{\tilde{M}}_{16}^{A}$
and $\mathbf{\tilde{M}}_{16}^{B}$ to pandiagonal magic squares by the Planck
transformation \cite{PLAN,NORD0}.

\section{Conclusion}

Commuting magic square matrices can be formed by compounding two magic
square matrices with the all-ones matrix in a known manner. We verify that
the compounded matrices retain the regular (associative) and pandiagonality
properties of the original magic squares as noted by Eggermont \cite{EGGE}.
In the case where a single matrix is compounded with the all-ones matrix in
two ways, the compounded matrices are related by a row/column permutation
(shuffle) of their elements as noted by Rogers, et.\negthinspace\ al \cite%
{ROGE} and expressed here in a matrix form which shows that they are
similar. We derive simple formulas for the Jordan-form matrices and SVD
matrices of the compounded magic square matrices in terms of those of the
original magic squares. A wider class of commuting magic squares is
considered but our formulas for the Jordan form and SVD do not apply to
them. Three examples illustrate the constructions and validate our formulas.
The methods presented here should apply to other compound matrix
constructions, such as the additional ones given in \cite{ROGE}.\smallskip 

\noindent \textbf{Acknowledgement.} I am grateful to Peter Loly for a
helpful discussion and for providing pertinent pages of Eggermont's thesis 
\cite{EGGE} which is no longer available on the internet or elsewhere.

\end{document}